\documentclass{article}
\usepackage{amsmath}
\usepackage{amsfonts}
\usepackage{xy}
\input xy
\xyoption{all}
\newcommand{\Supp}{{\mathrm {Supp}}}
\newcommand{\supp}{{\mathrm {supp}}}

\newcommand{\Star}{{\mathrm {Star}}}
\newcommand{\Pen}{{\mathrm {Pen}}}

\newcommand{\Cliff}{{\mathrm {Cliff}}}

\newcommand{\cA}{{\mathcal{A} }}
\newcommand{\cB}{{\mathcal{B} }}
\newcommand{\cC}{{\mathcal{C} }}

\newcommand{\cK}{{\mathcal{K} }}

\newcommand{\cM}{{\mathcal{M} }}

\newcommand{\cS}{{\mathcal{S} }}

\newcommand{\wt}{\widetilde}
\newcommand{\nn}{{n\in \mathbb{N}}}
\newcommand{\ttt}{T=\big[\big( T^{(0)}, \cdots, T^{(n)}, \cdots \big) \big]}

\newcommand{\pdx}{P_d(X)}
\newcommand{\pdxn}{P_d(X_n)}
\newcommand{\zdn}{{Z_{d,n}}}

\newcommand{\cx}{\mathbb{C}[X]}
\newcommand{\cmx}{C^*_{\max}(X)}
\newcommand{\cmy}{C^*_{\max}(Y)}

\newcommand{\cpdx}{\mathbb{C}\big[ P_d(X) \big] }
\newcommand{\cmpdx}{C^*_{\max}\big( \pdx \big) }

\newcommand{\cmpdxn}{C^*_{\max}\big( \pdxn \big) }

\newcommand{\cupdxn}{\mathbb{C}_{u, \infty} \big[ \big( P_d(X_n)\big)_{n\in \mathbb{N}} \big] }
\newcommand{\cumpdxn}{C^*_{u,\max,\infty} \big( \big(P_d(X_n)\big)_{n\in \mathbb{N}} \big) }

\newcommand{\culpdxn}{\mathbb{C}_{u,L, \infty} \big[ \big(P_d(X_n)\big)_{n\in \mathbb{N}} \big] }
\newcommand{\culmpdxn}{C^*_{u,L, \max,\infty} \big( \big(P_d(X_n)\big)_{n\in \mathbb{N}} \big) }

\newcommand{\cupdxna}{\mathbb{C}_{u, \infty} \big[ \big(P_d(X_n), \cA(V_n) \big)_{n\in \mathbb{N}} \big] }
\newcommand{\cumpdxna}{C^*_{u,\max,\infty} \big( \big(P_d(X_n), \cA(V_n) \big)_{n\in \mathbb{N}} \big) }

\newcommand{\culpdxna}{\mathbb{C}_{u,L,\infty} \big[ \big(P_d(X_n), \cA(V_n) \big)_{n\in \mathbb{N}} \big] }
\newcommand{\culmpdxna}{C^*_{u,L, \max,\infty} \big( \big( P_d(X_n), \cA(V_n) \big)_{n\in \mathbb{N}} \big) }

\newcommand{\cupdxnk}{\mathbb{C}_{u, \infty} \big[ \big(P_d(X_n), \cK(L^2_n) \big)_{n\in \mathbb{N}} \big] }
\newcommand{\cumpdxnk}{C^*_{u,\max,\infty} \big( \big(P_d(X_n), \cK(L^2_n) \big)_{n\in \mathbb{N}} \big) }

\newcommand{\culmpdxnk}{C^*_{u,L, \max,\infty} \big(\big( P_d(X_n), \cK(L^2_n) \big)_{n\in \mathbb{N}} \big) }
\newcommand{\ggnn}{\gamma\in \Gamma_n, \nn}
\begin{document}
	\title{The maximal coarse Baum-Connes conjecture for spaces which admit a fibred coarse embedding into Hilbert space
			\thanks{The first author is supported by NSFC (No. 10731020). The second author is supported by NSFC (No. 10971023),
					the Fundamental Research Funds for the Central Universities of China, and Shanghai Shuguang Project (No. 07SG38).
					The third author is supported by the U.S. National Science Foundation. }
		 }
	\author{Xiaoman Chen, Qin Wang and Guoliang Yu}
	\maketitle
\begin{abstract}
We introduce a notion of fibred coarse embedding into Hilbert space for metric spaces, which is a generalization of Gromov's notion of coarse embedding into Hilbert space. It turns out that a large class of expander graphs admit such an embedding. We show that the maximal coarse Baum-Connes conjecture holds for  metric spaces with bounded geometry which admit a fibred coarse embedding into Hilbert space.
\end{abstract}
\tableofcontents
\section{Introduction}
In this paper, we shall introduce a concept of {\em fibred coarse embedding into Hilbert space}  for metric spaces (Definition 2.1) and prove the following result:
\par
{\bf Theorem 1.1.} {\it Let $X$ be a discrete metric space with bounded geometry. If $X$ admits a fibred coarse embedding into Hilbert space, then the maximal coarse
Baum-Connes conjecture holds for $X$.}
\par
The (maximal) coarse Baum-Connes conjecture \cite{HR, Yu95, GWY} is a geometric analogue of the Baum-Connes conjecture \cite{BC, BCH} and provides an algorithm  for computing the higher
indices of generalized elliptic operators on non-compact spaces. It states
that a certain assembly map $\mu$ or $\mu_{\max}$ from $\lim_{d\to \infty} K_*(P_d(X))$ to $K_*(C^*(X))$ or $K_*(C^*_{\max}(X))$ is an isomorphism, where $K_*(P_d(X))$ is the
locally finite $K$-homology group
of the Rips complex $P_d(X)$,  and  $K_*(C^*(X))$ and $K_*(C^*_{\max}(X))$ are respectively the $K$-theory groups of the Roe algebra and the maximal Roe algebra of $X$.
The conjecture has many applications in topology and geometry. In particular, it implies the Novikov conjecture on homotopy invariance of higher signatures,
the Gromov conjecture on non-existence of positive scalar curvature metrics on uniformly contractible Riemannian manifolds,
and the zero-in-the-spectrum conjecture stating that the Laplacian operator acting on
the space of all $L^2$-forms of a uniformly contractible Riemannian manifold has zero in its spectrum.
\par
Recall that a discrete metric space $X$ is said to have {\em bounded geometry} if for any $r>0$ there is $N>0$ such that any ball of radius $r$ in $X$ contains at most $N$ elements.
A map $f: X\to Y$ from a metric space $X$ to another $Y$ is said to be a {\em coarse embedding} or {\em uniform embedding} \cite{Grom93} if there exist non-decreasing functions
$\rho_1$ and $\rho_2$ from $[0,\infty)$ to $(-\infty, \infty)$ with $\displaystyle\lim_{r\to \infty}\rho_i(r)=\infty$ ($i=1,2$) such that
$$\rho_1\big( d(x,y) \big) \leq d\big( f(x),f(y) \big)\leq \rho_2\big(d(x,y)\big)$$
for all $x, y\in X$.  \quad M. Gromov suggested that coarse embeddability of a discrete group into a Hilbert space, or a certain Banach space, might be relevant to prove the Novikov conjecture
\cite{Grom93}. G. Yu subsequently proved the coarse Baum-Connes conjecture for discrete metric spaces with bounded geometry which are coarsely embeddable into Hilbert space \cite{Yu00}.
M. Gromov discovered that  a metric space of expander graphs \cite{Marg, Lub} cannot be coarsely embedded into a Hilbert space \cite{Grom00}. N. Higson showed that the assembly map $\mu$ fails to be
surjective for certain Margulis-type of expanders \cite{Hig99}.
%
%
In \cite{Laf}, V. Lafforgue showed that there are residually finite groups whose
associated expander graphs cannot be coarsely embedded into any uniformly convex Banach space. Meanwhile, G. Gong, Q. Wang and G. Yu introduced the maximal Roe algebra
of a metric space with bounded geometry  in \cite{GWY}, and proved  a version of the maximal coarse Novikov conjecture, i.e. the injectivity of the maximal assembly map $\mu_{\max}$,
for the box space of a class of residually finite groups, including the expander graphs
constructed by V. Lafforgue.
Thereafter, the coarse Novikov conjecture, i.e. the injectivity of the assembly map $\mu$,
 for a large class of expanders is proved in \cite{CTWY, GTY, OY}.
More recently, H. Oyono-Oyono and G. Yu proved the maximal coarse Baum-Connes conjecture for certain metric spaces constructed from spaces with isometric actions by
residually finite groups \cite{OY}, including a certain class of expander graphs. R. Willett and G. Yu proved the maximal coarse Baum-Connes conjecture for spaces of graphs with large
girth \cite{WiY2}, also including a certain class of expander graphs.
\par
The notion of fibred coarse embedding into Hilbert space, which we are going to introduce,  generalizes Gromov's notion of coarse embedding into Hilbert space to a great extent.
Roughly speaking, a metric space $X$ admits such an embedding implies that, although the whole space $X$ may not be coarsely embedded into Hilbert space globally, large bounded subsets of  $X$
can be coarsely embedded into Hilbert space within a common distortion as long as these large bounded subsets are far away towards infinity.
This feature provides the property with much flexibility, allowing many expander graphs studied above to admit such an embedding.
In particular, our result Theorem 1.1 includes the major results of \cite{OY} and \cite{WiY2} as special cases.
\par
This paper is organized as follows. In section 2, we introduce the notion of fibred coarse embedding into Hilbert space, and discuss sevaral situations to which this notion applies.
In section 3, we recall the definition of maximal Roe algebra and the formulation of the maximal coarse Baum-Connes conjecture. In section 4, we explain the strategy to prove Theorem 1.1.
The problem of proving the maximal coarse Baum-Connes conjecture for a metric space $X$ can be reduced to verifying that the evaluation homomorphism from the $K$-theory of a certain localization
algebra to the $K$-theory of the maximal Roe algebra at infinity for the coarse disjoint union of a sequence of finite subspaces of $X$ is an isomorphism.
In section 5, we define the twisted Roe algebra at infinity and its localization algebra for a sequence of finite metric spaces
which admits a fibred coarse embedding into Hilbert space. In section 6, we discuss various ideals of the twisted algebras and show that the evaluation map for the twisted algebras
is an isomorphism. In section 7, we prove a geometric analogue of the Bott periodicity in finite dimension by constructing the Bott maps and the Dirac maps as asymptotic morphisms. This is used to show
the evaluation map required in section 4 is an isomorphism, which implies Theorem 1.1.
\par
Throughout this paper, we denote by $\cK(H)$ or $\cB(H)$ the algebras of all compact operators or all bounded linear operators on a Hilbert space $H$. We denote by
$\cK$ the algebra of all compact operators on a fixed separable Hilbert space $H_0$. Denote by $\mathbb{N}=\{0, 1, 2, \cdots\}$ the non-negative integers. Denote by $\#A$ the number of elements in a set $A$.
For a metric space $X$, a point $x\in X$ and $r>0$, we denote the ball of  $X$ with center $x$ and radius $r$  by $B(x, r)$,  or $B_X(x, r)$ if it is necessary to indicate the base space $X$.

\par
{\bf Acknowledgement.} The authors are very grateful to Rufus Willett who carefully read an early version of this paper and suggested very helpful comments.

\section{Fibred coarse embedding into Hilbert space}
In this section, we introduce the concept of {\em fibred coarse embedding into Hilbert space} for a metric space, generalizing Gromov's notion of
{\em coarse embedding into Hilbert space} \cite{Grom93}. Let $H$ be a separable Hilbert space, as a model space.
\par
{\bf Definition 2.1.} A metric space $(X, d)$ is said to admit a {\em fibred coarse embedding into Hilbert space} if there exist
\begin{itemize}
\item a field of Hilbert spaces $(H_x)_{x\in X}$ over $X$;
\item a section $s: X\to \bigsqcup_{x\in X} H_x$ (i.e. $s(x)\in H_x$);
\item two non-decreasing functions $\rho_1$ and $\rho_2$ from $[0, \infty)$ to $(-\infty, \infty)$ with $\lim_{r\to \infty}\rho_i(r)=\infty$ ($i=1, 2$)
\end{itemize}
such that for any $r>0$ there exists a bounded subset $K\subset X$ for which there exists a ``trivialization''
$$t_C: (H_x)_{x\in C}\longrightarrow C\times H$$
for each subset $C\subset X\backslash K $ of diameter less than $r$, i.e. a map from $(H_x)_{x\in C}$ to the constant field $C\times H$ over $C$ such that the restriction  of $t_C$ to
the fiber $H_x$ ($x\in C$) is an affine isometry $t_C(x): H_x\to H$, satisfying
\begin{itemize}
\item[(1)] for any $x, y\in C$, $\rho_1(d(x, y))\leq \|t_C(x)(s(x))-t_C(y)(s(y))\|\leq \rho_2(d(x, y))$;
\item[(2)] for any two subsets $C_1, C_2\subset X\backslash K$ of diameter less than $r$ with $C_1\cap C_2\not=\emptyset$, there exists an affine isometry $t_{C_1C_2}: H\to H$ such that
		$t_{C_1}(x)\circ t_{C_2}^{-1}(x)=t_{C_1C_2}$ for all $x\in C_1\cap C_2$.
\end{itemize}
\hfill{$\Box$}
\par
In the following we discuss several situations to which the above notion applies. Let $(X_n)_\nn$ be a sequence of bounded metric spaces.
\par
Recall that {\em a coarse disjoint union of $(X_n)_\nn$} is the disjoint union $X=\bigsqcup_\nn X_n$ equipped with a metric $d$ such that (1) the restriction of $d$ to each $X_n$ is the original metric
of $X_n$; (2) $d(X_n, X_m)\to \infty$ as $n+m\to \infty$ and $n\neq m$. Note that any two metrics satisfying these two conditions are coarsely equivalent. So we may refer to $X$ as {\em the} coarse disjoint
union of  $(X_n)_\nn$.
\par
For each $\nn$, a metric space $\wt{X_n}$ is called {\em a Galois covering of $X_n$} if there exists a discrete group $\Gamma_n$ acting on $\wt{X_n}$ freely and properly by isometries
such that $X_n=\wt{X_n}/\Gamma_n$. Denote by $\pi_n: \wt{X_n}\to X_n$ the associated covering map.
\par
A sequence of Galois coverings $(\wt{X_n})_\nn$ of $(X_n)_\nn$ is said to be {\em asymptotically faithful} (cf. \cite{WiY1}) if for any $r>0$ there exists $N\in \mathbb{N}$ such that, for all
$n\geq N$, the covering map $\pi_n: \wt{X_n}\to X_n$ is ``$r$-isometric'', i.e. for all subsets $\wt{C}\subset \wt{X_n}$ of diameter less than $r$, the restriction of $\pi_n$ to $\wt{C}$ is an isometry
onto $C=\pi_n(\wt{C})\subset X_n$.
\par
A sequence of Galois coverings $(\wt{X_n})_\nn$ of $(X_n)_\nn$ is said to admit a {\em uniform equivariant coarse embedding into Hilbert space} if there exist a map $f_n: \wt{X_n}\to H$ for each
$\nn$, and two non-decreasing functions $\rho_1$ and $\rho_2$ from $[0, \infty)$ to $(-\infty, \infty)$ with $\lim_{r\to \infty}\rho_i(r)=\infty$ ($i=1, 2$), such that
(1) $f_n$ is $\Gamma_n$-equivariant (this implies $H$ has an action of $\Gamma_n$ by affine isometries) for all $\nn$; (2) $\rho_1(d(x, y))\leq \|f_n(x)-f_n(y)\|\leq \rho_2(d(x, y))$
for all $x, y\in \wt{X_n}$, $\nn$.
\par
{\bf Theorem 2.2.} {\it
Let $X=\bigsqcup_\nn X_n$ be a coarse disjoint union of a sequence of bounded metric spaces. If there exists a sequence of asymptotically faithful Galois coverings
$(\wt{X_n})_\nn$ which admits a uniform equivariant coarse embedding into Hilbert space, then $X$ admits a fibred coarse embedding into Hilbert space.
}
\par
{\bf Proof.} For each $\nn$,  the group $\Gamma_n$ acts on $\wt{X_n}\times H$ by
$$\gamma \cdot (x, h)= (\gamma x, \gamma h)$$
for $x\in \wt{X_n}, h\in H, \gamma\in \Gamma_n$. Denote the orbit of $(x, h)\in \wt{X_n}\times H$ by $[(x, h)]$, i.e.
$$[(x, h)]=\{ \, (\gamma x, \gamma h) \; | \; \gamma\in \Gamma_n \, \}.$$
For any point $a\in X_n=\wt{X_n}/\Gamma_n$, the action
of a group element $\gamma\in \Gamma_n$ on $\wt{X_n}$ permutes the points in $\pi_n^{-1}(a)\subset \wt{X_n}$. Define
$$H_a:=\Big( \pi_n^{-1}(a) \times H \Big)\Big/ \Gamma_n.$$
Then
$$ (H_a)_{a\in X_n}=\Big( \wt{X_n} \times H \Big)\Big/\Gamma_n $$
is a field of Hilbert spaces over $X_n$. Define a section
$$s: X_n\longrightarrow \Big( \wt{X_n} \times H \Big)\Big/\Gamma_n$$
by the formula
$$s(a)=[(x, f_n(x))]\in H_a$$
for any $x\in \pi_n^{-1}(a)$, $a\in X_n$. This is well-defined since the map $f_n: \wt{X_n}\to H$ is $\Gamma_n$-equivariant.
\par
Since the covering sequence $(\wt{X_n})_\nn$ is asymptotically faithful, for any $r>0$ there exists $N\in \mathbb{N}$ such that
for any $n\geq N$ and any $C\subset X_n$ of diameter less than $r$,
the action of an element $\gamma \in \Gamma_n$ on $\wt{X_n}$ permutes the (disjoint) components of $\pi_n^{-1}(C)$, and the restriction of $\pi_n$ to each of these components is an
isometry onto $C$. For a point $z\in \pi_n^{-1}(C)$, denote by ${\wt{C}}^z$ the component of $\pi_n^{-1}(C)$ containing $z$. Then each such
component ${\wt{C}}^z$ gives rise to a trivialization:
$$t_{C, z}: (H_a)_{a\in C}=\Big( \pi_n^{-1}(C) \times H \Big)\Big/ \Gamma_n  \longrightarrow C\times H, $$
where, for any $a\in C\subset X_n$, the affine isometry
$$t_{C, z}(a): \; \; H_a=\Big( \pi_n^{-1}(a) \times H \Big)\Big/ \Gamma_n     \longrightarrow H$$
is given by
$$\Big( t_{C, z}(a) \Big) \Big( [(x, h)] \Big) =h$$
for $[(x, h)]\in H_a$ represented by $(x, h)\in \wt{C}^z\times H$.
\par
Now, for any $a, b\in C$, there exist $x, y\in \wt{C}^z$ such that $\pi_n(x)=a$ and $\pi_n(y)=b$, so that
\[
\begin{array}{c}
  \Big\|\Big(t_{C, z}(a)\Big) \big( s(a) \big) - \Big( t_{C, z}(b) \Big) \big( s(b) \big) \Big\| = \Big\| f_n(x)-f_n(y) \Big\|              \\
\in  \big[ \rho_1(d(x, y)), \; \rho_2(d(x, y)) \big] =  \big[ \rho_1(d(a, b)), \; \rho_2(d(a, b)) \big].
\end{array}
\]
\par
Moreover, for any $C_1, C_2\subset X_n$ ($n\geq N$) of diameter less than $r$ with $C_1\cap C_2\not=\emptyset$, there exist $z, \gamma z\in \wt{X_n}$ for some $\gamma\in \Gamma_n$ such that
$\pi_n(z)=\pi_n(\gamma z)\in C_1\cap C_2$. Then,  for all $a\in C_1\cap C_2$,
$$t_{C_1, z}(a)\circ t_{C_2, \gamma z}^{-1}(a)=\gamma: \; H\to H$$
as an affine isometry mapping $h\in H$ to $\gamma h\in H$. The proof is complete.
\hfill{$\Box$}
\par
{\bf Example 2.3.} If a  metric space $X$ is coarsely embeddable into Hilbert space, then clearly it admits a fibred coarse embedding into Hilbert space. The maximal coarse Baum-Connes conjecture in this
case follows from G. Yu's  work on the coarse Baum-Connes conjecture \cite{Yu00}, together with a recent result by J. \v{S}pakula and R. Willett that $K_*(C^*_{\max}(X))$ is isomorphic to $K_*(C^*(X))$ if $X$ is coarsely embeddable into Hilbert space \cite{SpWi}.  Note also that for a discrete group, fibred coarse embeddability is the same as usual coarse embeddability.
\hfill{$\Box$}
\par
{\bf Example 2.4.} Let $X$ be a discrete metric space with bounded geometry. Suppose a residually finite group $\Gamma$ acts on $X$ freely, properly and cocompactly by isometries,
 such that $X$ is $\Gamma$-equivariantly coarsely embeddable into Hilbert space. Let $(\Gamma_n)_\nn$ be a sequence of finite index normal subgroups of $\Gamma$ such that for all $r>0$ there
exists $N\in \mathbb{N}$ so that if $B(e, r)$ is the ball in $\Gamma$ of radius $r$ about the identity, then $\Gamma_n\cap B(e, r)=\{e\}$ for all $n\geq N$.
Let $X_n=X/\Gamma_n$ be the quotient space. It follows from Theorem 2.2 that the coarse disjoint union $X=\bigsqcup_\nn X_n$ admits a fibred coarse embedding into Hilbert space. Indeed, the constant sequence
$\wt{X_n}=X$ ($\nn$)
serves as the asymptotically faithful Galois coverings.  Theorem 1.1 implies the maximal coarse Baum-Connes conjecture for such X.
A special case of this result was proved by H. Oyono-Oyono and G. Yu in \cite{OY}.
\par
For a residually finite group $\Gamma$, the space $X(\Gamma)=\bigsqcup_{\nn} \Gamma/\Gamma_n$ is called {\em the box space of  $\Gamma$} (\cite{Roe03}). It turns out that if the box space $X(\Gamma)$ is coarsely embeddable into Hilbert space, then $\Gamma$ is a-T-menable
(\cite{Roe03}). It is also easy to see that the converse of this implication is not true. However, it is shown in \cite{CWW12} that $\Gamma$ is a-T-menable if and only if the box space $X(\Gamma)$ is fibred coarsely embeddable into Hilbert space.
\hfill{$\Box$}
\par
{\bf Example 2.5.} Recall that the {\em girth} of a graph $G$, denoted by $girth(G)$,  is the shortest length of a cycle in $G$. A sequence of finite connected graphs $(G_n)_\nn$ is
said to have {\em large girth} if $girth(G_n)\to \infty$ as $n\to \infty$. Then the coarse disjoint union $X=\bigsqcup_{\nn} G_n$ admits a fibred coarse embedding into Hilbert space.
Indeed, let $\wt{G_n}$ be the {\em universal cover of $G_n$}, which is actually a tree. Then the covering sequence $(\wt{G_n})_\nn$ satisfies the conditions in Theorem 2.2. The maximal coarse Baum-Connes
conjecture for such $X$ was proved by R. Willett and G. Yu in \cite{WiY2}.
\hfill{$\Box$}
\par
Note that Examples 2.4 and 2.5 imply that a large class of expander graphs admit a fibred coarse embedding into Hilbert space.

\section{The maximal coarse Baum-Connes conjecture}

In this section, we shall collect from \cite{GWY} results concerning the maximal Roe algebra of a proper metric space with bounded geometry and the formulation of the maximal coarse Baum-Connes conjecture
(see also \cite{OY, WiY1, WiY2}).
\par
Let $X$ be a proper metric space (a metric space is called $proper$ if every closed ball is compact). An $X$-module $H_X$ is a separable Hilbert space
equipped with a  $*$-representation $\pi$ of $C_0(X)$, the algebra of all continuous functions on $X$ which vanish
at infinity. An $X$-module is called non-degenerate if the $*$-representation of $C_0(X)$ is non-degenerate.
An $X$-module is said to be standard if no nonzero function in $C_0(X)$ acts as a compact operator.
When $H_X$ is an $X$-module,  for each $f\in C_0(X)$ and $h\in H_X$, we denote $(\pi (f))h$ by $fh$.

\par
{\bf Definition 3.1.} (cf. \cite{Roe93})  Let $X$ be a standard non-degenerate $X$-module.
\\ \indent (1) The {\em support} $\supp(T)$ of a bounded linear operator $T: H_X\to H_X$ is defined to be the complement of the set of all points $(x, y)\in X\times X$ for which
there exist $f, g\in C_0(X)$  such that $gTf=0$ but $f(x)\not= 0$, $g(y)\not=0$.
\\ \indent (2) A bounded operator $T: H_X\to H_X$  is said to have {\em finite propagation}  if
$$\sup\{d(x, y): (x, y)\in \Supp(T)\}<\infty.$$
This number is called the {\em propagation of $T$}.
\\ \indent (3) A bounded operator $T: H_X\to H_X$ is said to be {\em locally compact} if the operators $fT$ and $Tf$ are compact for all $f\in C_0(X)$.
\par
Denote by $\mathbb{C}[X, H_X]$, or simply $\mathbb{C}[X]$, the set of all locally compact, finite propagation operators on a standard non-degenerate $X$-module $H_X$.  It is straightforward to check that $\mathbb{C}[X]$ is a $*$-algebra which,
up to non-canonical isomorphisms, does not depend on the choice of standard non-degenerate $X$-module.
\par
{\bf Definition 3.2.} A {\em net} of a metric space $X$ is a countable subset $\Gamma\subset X$ such that there exist numbers $\delta>0$, $R>0$ satisfying (1) $d(\gamma, \gamma')>\delta$
for all distinct elements $\gamma, \gamma'\in \Gamma$; (2) for any $x\in X$ there exists $\gamma\in \Gamma$ such that $d(x, \gamma)<R$.  \quad A metric space $X$ is said to have {\em bounded geometry}
if $X$ contains a net with bounded geometry.
\par
The following result is essentially proved in \cite{GWY}.
\par
{\bf Lemma 3.3.} {\it Let $X$ be a proper metric space with bounded geometry, and let $H_X$ be a standard non-degenerate $X$-module. For any $r>0$ there exists a constant $c>0$ such that
for any $*$-representation $\phi$ of $\cx$ on a Hilbert space $H_\phi$ and any $T\in \cx$ with propagation less than $r$, we have
$$||\phi(T)||_{\cB(H_\phi)}\leq c  \, || T ||_{\cB(H_X)} \,. $$
}
\par
This allows us to define the maximal Roe algebra of $X$.
\par
{\bf Definition 3.4.} \quad (\cite{GWY}) Let $X$ be a proper metric space with bounded geometry. The {\em maximal Roe algebra} of $X$, denoted by $\cmx$, is the completion of
$\cx$ with respect to the $C^*$-norm:
$$\big\|T\big\|_{\max} :=\sup\Big\{ \big\|\phi(T)\big\|_{\cB(H_\phi)} \; \Big | \; \phi:\cx\to \cB(H_\phi), \mbox{ a $*$-representation} \Big\}. $$
\par
{\bf Definition 3.5.}
Let $H_X=\ell^2(Z, H_0)$, where $Z\subset X$ is a countable dense subset of $X$ and $H_0$ is an infinite dimensional separable Hilbert space. A function $f\in C_0(X)$ acts on
$\ell^2(Z)\otimes H_0$ by pointwise multiplication on the first component: $f(\xi\otimes h)=f\xi\otimes h$. Define $\mathbb{C}_f [X]$ to be the $*$-algebra of all bounded functions
$T: Z\times Z\to \cK:=\cK(H_0)$, also viewed as $Z\times Z$-matrices, such that
\begin{itemize}
\item for any bounded subset $B\subset X$, the set
		$$\{(x, y)\in B\times B\cap Z\times Z \; | \; T(x, y)\not= 0 \}$$
		 is finite;
\item there exists $L>0$ such that
		$$\#\{y\in Z \; | \; T(x, y)\not= 0 \}<L,  \quad \#\{y\in Z \; | \; T(y, x)\not= 0 \}<L $$
		for all $x\in Z$;
\item there exists $R>0$ such that $T(x, y)=0$ whenever $d(x, y)>R$ for $x, y\in Z$.
\end{itemize}
Note that in general $\mathbb{C}_f [X]$ is a dense $*$-subalgebra of $\mathbb{C}[X, \ell(Z, H_0)]$ within $C^*_{\max}(X)$. In the sequel, we will use $\mathbb{C}_f [X]$ to replace $\mathbb{C}[X]$ as a generating subalgebra of $C^*_{\max}(X)$.
\par
{\bf Remark 3.6.} (cf. \cite{HRY} and \cite{GWY}) \quad (1) $\cmx$ does not depend on the choice of countable dense subset $Z\subset X$ up to non-canonical isomorphisms.
If $X$ and $Y$ are coarsely equivalent, then $\cmx$ is isomorphic to $\cmy$ via a non-canonical isomorphism. (2) The $K$-theory groups of $\cmx$ do not depend on the choice of $Z$ up to
canonical isomorphisms. If $X$ and $Y$  are coarsely equivalent, then $K_*(\cmx)$ is isomorphic to $K_*(\cmy)$ via a canonical isomorphism.
\par
We next define the assembly map $\mu_{\max}$ (also referred to as the index map or the Baum-Connes map) for the maximal Roe algebras.
\par
Let $X$ be a proper metric space. Recall that the locally finite $K$-homology groups $K_i(X)=KK_i(C_0(X),\mathbb{C})$ $(i=0,1)$ are generated by certain cycles
(abstract elliptic operators) modulo  certain equivalence relations \cite{Kas75, Kas88}:
\\ \indent (1) a cycle for $K_0(X)$ is a pair $(H_X, F)$, where $H_X$ is an $X$-module and $F$ is a bounded linear operator acting on $H_X$ such that $F^*F-I$ and $FF^*-I$ are locally compact,
and $\phi F-F\phi$ is compact for all $\phi\in C_0(X)$;
\\ \indent (2) a cycle for $K_1(X)$ is a pair $(H_X, F)$, where $H_X$ is an $X$-module and $F$ is a {\em self-adjoint} operator acting on $H_X$ such that $F^2-I$ is locally compact,
and $\phi F-F\phi$ is compact for all $\phi\in C_0(X)$.
\par
In both cases, the equivalence relations on cycles are given by homotopy of the operators $F$, unitary equivalence, and direct sum with ``degenerate" cycles, those cycles for which
$F\phi-\phi F$, $\phi(F^* F-I)$ and so on, are not merely compact but actually zero \cite{Kas75, Kas88}.
\par
The assembly map $\mu_{\max}: K_*(X) \to K_*(\cmx)$ is defined as follows.
\par
{\bf Definition 3.7.} Let $(H_X, T)$ represent a cycle in $K_0(X)$. For any $R>0$, one can always take a  locally finite,
uniformly bounded open cover $\{U_i\}_i$ of $X$ such that
$$diameter(U_i)<R $$
for all $i$, and a continuous partition of unity  $\{\phi_i\}_i$ subordinate to the open cover $\{U_i\}_i$. Define $F=\sum_i \phi_i^{\frac{1}{2}}T\phi_i^{\frac{1}{2}}$, where the sum converges in
the strong operator topology. It is not hard to see that $(H_X, T)$ and $(H_X, F)$ are equivalent in $K_0(X)$. Note that now the propagation of $F$ is less than or equal to $R$, so that
$F^*F-I$ and $FF^*-I$ are in $\mathbb{C}[X]$. Let (cf. \cite{Mil})
\[W=
\left(
\begin{array}{cc}
I & F \\
0 & I
\end{array}
\right)
\left(
\begin{array}{cc}
I & 0 \\
-F^* & I
\end{array}
\right)
\left(
\begin{array}{cc}
I & F \\
0 & I
\end{array}
\right)
\left(
\begin{array}{cc}
0 & -I \\
I & 0
\end{array}
\right)
\in \cB(H_X\oplus H_X).
\]
Then both $W$ and $W^{-1}$ have finite propagation (at most $3R$), and
\[
W
\left(
\begin{array}{cc}
I & 0 \\
0 & 0
\end{array}
\right)
W^{-1}-
\left(
\begin{array}{cc}
I & 0 \\
0 & 0
\end{array}
\right)
\in \mathbb{C}[X]\otimes \cM_2(\mathbb{C}).
\]
We then define
\[
\mu \big( [(H_X, T)] \big):=
\left[
W
\left(
\begin{array}{cc}
I & 0 \\
0 & 0
\end{array}
\right)
W^{-1}
\right]
-
\left[
\left(
\begin{array}{cc}
I & 0 \\
0 & 0
\end{array}
\right)
\right]
\]
in $K_0(\mathbb{C}[X])$. Furthermore, $\mu \big( [(H_X, T)] \big)$ defines an element in $K_0(C^*_{\max}(X))$ by considering $\mathbb{C}[X]$ as a  $*$-subalgebra of $C^*_{\max}(X)$. This element is denoted by $\mu_{\max} \big( [(H_X, T)] \big) \in K_0(C^*_{\max}(X))$.
Thus, we obtain the assembly map $\mu_{\max}: K_0(X)\rightarrow K_0(C^*_{\max}(X))$. Similarly, we can define $\mu_{\max}: K_1(X)\rightarrow K_1(C^*_{\max}(X))$.
\par
{\bf Remark 3.8.}
Let $X$ be a proper metric space with bounded geometry as above, and $Z$ a countable dense subset of $X$ used to define $\mathbb{C}_f [X]$ as in Definition 3.5. For any natural number $n>0$, let $Z_n$ be a subset of $Z$ such that $d(x, y)>\frac{1}{2n}$ for distinct $x, y\in Z_n$, and $d(x, Z_n)\leq \frac{1}{n}$ for all $x\in X$. Without loss of generality, we may assume $Z=\cup_{n=1}^\infty Z_n$. Let $\mathbb{C}[Z_n]$ be the $*$-algebra of all bounded functions $T: Z_n\times Z_n\to \mathcal{K}$ with finite propagation, i.e., there exists $R>0$ such that $T(x, y)=0$ whenever $d(x, y)>R$ for all $x, y\in Z_n$. Then $\mathbb{C}[Z_n]\subset \mathbb{C}_f [X]$. Moreover, there exists a non-canonical $*$-isomorphism (cf. \cite{HRY} or 4.4 in \cite{GWY}) $Ad(U): \mathbb{C}[X]\to \mathbb{C}[Z_n]$ such that, if $T\in \mathbb{C}[X]$ has propagation less than $R$, then the propagation of $(Ad(U))(T)$ is less than $R+\frac{2}{n}$.
\par
For any $R>0$, let $\mu \big( [(H_X, T)] \big)\in K_0(\mathbb{C}[X])$ be as in Definition 3.7. Then
$$Ad(U)_*\big(\mu \big( [(H_X, T)] \big) \big)\in K_0(\mathbb{C}[Z_n]),$$
which also defines an element in $K_0(\mathbb{C}_f[X])$ via the inclusion $\mathbb{C}[Z_n]\hookrightarrow \mathbb{C}_f [X]$.
Note that the propagation of the element
\[
\big( Ad(U) \big)\left(
W
\left(
\begin{array}{cc}
I & 0 \\
0 & 0
\end{array}
\right)
W^{-1}-
\left(
\begin{array}{cc}
I & 0 \\
0 & 0
\end{array}
\right)\right)
\in \mathbb{C}_f[X]\otimes \cM_2(\mathbb{C})
\]
is less than $6R+\frac{2}{n}$. Note also that $Ad(U)_*$ canonically induces the identity on $K_*(C^*_{\max}(X))$.
\par
{\bf Definition 3.9.} Let $X$ be a discrete metric space with bounded geometry. For each $d\geq 0$, the {\em Rips complex} $P_d(X)$ at scale $d$ is defined to be the simplicial polyhedron
in which the set of  vertices is $X$, and a finite subset $\{x_0, x_1, \cdots, x_q\}\subset X$ spans a simplex if and only if $d(x_i,x_j)\leq d$ for all $0\leq i, j\leq q$.
\par
Endow $P_d(X)$ with the {\em spherical metric}. Recall that on each path connected component of $P_d(X)$, the spherical
metric is the maximal metric whose restriction to each simplex
$\{\sum_{i=0}^q t_ix_i| t_i\geq 0, \sum_{i=0}^q t_i=1\}$ is the metric obtained by identifying the simplex with
$S_{+}^q$ via the map
$$\sum_{i=0}^q t_ix_i \mapsto \left( \frac{t_0}{\sqrt{\sum_{i=0}^q t_i^2}}, \frac{t_1}{\sqrt{\sum_{i=0}^q t_i^2}}, \cdots,
\frac{t_q}{\sqrt{\sum_{i=0}^q t_i^2}} \right)
$$
where $S_+^q:=\{(s_0, s_1, \cdots, s_q)\in \mathbb{R}^{q+1}, s_i\geq 0, \sum_{i=0}^q s_i =1\}$ is endowed with the standard
Riemannian metric. If $y_0, y_1$ belong to two different connected components $Y_0, Y_1$ of $P_d(X)$, we define
$$d(y_0, y_1)= \min \{d(y_0, x_0)+d_X(x_0, x_1)+d(x_1, y_1)| x_0\in X\cap Y_0, x_1\in X\cap Y_1\}.$$
The topology induced by the above metric is the same as the weak topology of the simplicial complex:
a subset $S\subset P_d(X)$ is closed if and only if the intersection of $S$ with each simplex is closed.
\par
Note that for any $d\geq 0$, $\pdx$ is coarsely equivalent to $X$ via the inclusion map. If $d<d'$, then $P_d(X)$ is included in $P_{d'}(X)$ as a subcomplex via a simplicial map.
Passing to inductive limit, we obtain the assembly map
$$ \mu_{\max}: \lim_{d\to \infty} K_*(\pdx) \to  \lim_{d\to \infty} K_*(\cmpdx) \cong K_*(\cmx) .$$
\par
{\bf The maximal coarse Baum-Connes conjecture:}  {\it
If $X$ is a discrete metric space with bounded geometry, then the assembly map $\mu_{\max}$ is an isomorphism.
}
\par
{\bf Remark 3.10.} The maximal coarse Baum-Connes conjecture provides a method to compute the higher indices of elliptic differential operators
(Dirac operators) on non-compact complete Riemannian manifolds and has many applications in topology and geometry. In particular, it implies the Novikov conjecture, the Gromov positive scalar curvature conjecture,
etc. (cf. e.g. \cite{HR, Roe93, Roe96, Yu95, Yu98, Yu00, Yu06, GWY}).

\section{Reduction to coarse disjoint union of finite metric spaces}

The strategy to prove Theorem 1.1 is as follows. Let $X$ be a discrete metric space with bounded geometry which admits a fibred coarse embedding into Hilbert space.
We partition $X$ as $X=X^{(0)}\cup X^{(1)}$ such that each of $X^{(0)}$, $X^{(1)}$ and $X^{(0)}\cap X^{(1)}$ is a coarse disjoint union of a sequence of finite subspaces
of $X$, and the excision pair $(X^{(0)}, \; X^{(1)})$ respects the Mayer-Vietoris exact sequences on both the $K$-homology for these subspaces and the $K$-theory for the
corresponding maximal Roe algebras (the method of cutting the whole space into coarse disjoint unions of finite subspaces has its origin in \cite{WiY-book}). Consequently, it suffices to prove Theorem 1.1 for $X=\bigsqcup_{n\in \mathbb{N}} X_n$ being a coarse disjoint union of a sequence of finite
metric spaces. This case in turn can be reduced to verifying the evaluation homomorphism from the $K$-theory groups of certain localization algebras to that of the corresponding
{\em maximal Roe algebra at infinity for the sequence $(X_n)_\nn$} to be an isomorphism, a fact which will be proved in sections 5, 6 and 7. 
\par
In this section,  we shall explain in detail the above process of reduction. We start with the story for sequences of finite metric spaces. Let $(X_n)_\nn$ be a sequence of finite (discrete) metric spaces with {\em uniform bounded geometry},
i.e. for any $r>0$ there exists $N>0$ such that $\# B(x, r)<N$ for all $x\in X_n, \nn$. For each $d\geq 0$, let $P_d(X_n)$ be the Rips complex of $X_n$ at scale $d$ endowed with
the spherical metric. Take a countable dense subset $\zdn \subset \pdxn$ for each $d\geq 0$ in such a way that $\zdn\subseteq Z_{d', n}$ whenever $d<d'$, for all $\nn$.
\par
{\bf Definition 4.1.} For each $d\geq 0$, define $\cupdxn$ to be the set of all equivalence classes $\ttt$ of sequences $(T^{(0)}, \cdots, T^{(n)}, \cdots)$ described as follows:
\begin{itemize}
\item[(1)] $T^{(n)}$ is a bounded function from $\zdn \times \zdn$ to $\cK$ for all $\nn$;
\item[(2)] for any bounded subset $B\subset \pdxn$, the set
			$$\{ (x, y)\in B\times B \cap \zdn\times \zdn \; | \; T^{(n)} (x, y)\not= 0 \}$$
			is finite;
\item[(3)] there exists $L>0$ such that
			$$\#\{y\in \zdn \; | \; T^{(n)}(x, y)\not= 0 \}<L,  \quad  \#\{y\in \zdn \; | \; T^{(n)}(y, x)\not= 0 \}<L$$
			for all $x\in \zdn, \nn$;
\item[(4)] there exists $R>0$ such that $T^{(n)}(x, y)=0$ whenever $d(x, y)>R$ for $x, y\in \zdn, \; \nn$.
\end{itemize}
The equivalence relation $\sim$ on these sequences is defined by
$$(T^{(0)}, \cdots, T^{(n)}, \cdots) \sim (S^{(0)}, \cdots, S^{(n)}, \cdots)$$
 if and only if
$$\lim_{n\to \infty} \sup_{x, y\in \zdn} \big\| T^{(n)}(x, y) - S^{(n)}(x, y) \big\|_\cK =0.$$
Viewing $T^{(n)}$ as $\zdn\times \zdn$ matrices, $\cupdxn$ is then made into a $*$-algebra by using the usual matrix operations.
\par
Define $\cumpdxn$ to be the completion of $\cupdxn$ with respect to the norm
\begin{footnotesize}
$$\big\|T\big\|_{\max} :=\sup\Big\{ \big\|\phi(T)\big\|_{\cB(H_\phi)} \; \Big | \; \phi:\cupdxn \to \cB(H_\phi), \mbox{ a $*$-representation} \Big\}. $$
\end{footnotesize}
\hfill{$\Box$}
\par
Note that $\|T\|_{\max}$ is well-defined since $(X_n)_\nn$ have uniform bounded geometry. Moreover, $(\pdxn)_\nn$ is uniformly coarsely equivalent to $(X_n)_\nn$ for
any $d>0$, so that
$\cumpdxn$ is isomorphic to $C^*_{u, \max, \infty}((X_n)_\nn)$ via a non-canonical isomorphism. Recall from Definition 3.7 and Remark 3.8 that the individual assembly maps
$$\mu_{\max}: K_*(\pdxn) \to K_*(\cmpdxn)$$
($\nn$) can be defined by $\mu_{\max} ([(H_X, T)])=\mu_{\max} ([(H_X, F)])$ such that the propagation of $F$ is less than any given small $R>0$ which is independent of $\nn$, and that $\|F\|\leq 1$. Consequently,
we can obtain the following assembly map at infinity:
$$\mu_{\max, \infty}: \; \; \frac{\prod_{n=0}^\infty K_*(\pdxn)}{\oplus_{n=0}^\infty K_*(\pdxn)} \longrightarrow  K_*\Big(  \cumpdxn  \Big).$$
\par
The following notion of localization algebra has its origin in \cite{Yu97}.
\par
{\bf Definition 4.2.}  For each $d\geq 0$, define $\culpdxn$  to be the $*$-algebra of all bounded and uniformly norm-continuous functions
$$f: \; [0,\infty)\longrightarrow \cupdxn$$
such that $f(t)$ is of the form $f(t)=[(f^{(0)}(t), \cdots, f^{(n)}(t), \cdots )]$ for all $t\in [0, \infty)$, where the family of functions
$(f^{(n)}(t))_{\nn, t\geq 0}$ satisfy the conditions in Definition 4.1 with uniform constants, and there exists a bounded function $R:\; [0, \infty)\to [0,\infty)$ with
$\lim_{t\to \infty} R(t)=0$ such that
$$ \big( f^{(n)}(t)\big) (x, y) =0  \quad \mbox{whenever} \quad d(x, y)>R(t)$$
for all $x, y\in \zdn$, $\nn$, $t\in [0, \infty)$.
\par
Define $\culmpdxn$ to be the completion of $\culpdxn$ with respect to the norm
$$\big\| f \big\|_{\max} :=\sup_{t\in [0, \infty)} \big\|f(t)\big\|_{\max}.$$
\hfill{$\Box$}
\par
We can also define a local assembly map at infinity as in \cite{Yu97}:
$$\mu_{L, \max,\infty}: \; \; \frac{\prod_{n=0}^\infty K_*(\pdxn)}{\oplus_{n=0}^\infty K_*(\pdxn)} \longrightarrow  K_*\Big(  \culmpdxn  \Big).$$
\par
{\bf Proposition 4.3.}  {\it
Suppose $(X_n)_\nn$ have uniform bounded geometry. Then the local assembly map $\mu_{L, \max,\infty}$ is an isomorphism.
}
\hfill{$\Box$}
\par
The proof is similar to the arguments in \cite{Yu97}. Note that the evaluation homomorphism
$$e: \; \culmpdxn \longrightarrow \cumpdxn $$
defined by $e(f)=f(0)$ induces the following commutative diagram:

\[
\xymatrix{
 & \;\;\;\;\;\; \displaystyle    \lim_{d\to\infty}       K_*\big(   \culmpdxn  \big) \ar[d]^{e_*}            \\
\displaystyle\lim_{d\to\infty}     \frac{\prod_{n=0}^\infty K_*(\pdxn)}{\oplus_{n=0}^\infty K_*(\pdxn)}    \;\; \;\;\;  \ar[ur]^{\mu_{L, \max, \infty}}_\cong  \ar[r]^{\mu_{\max, \infty}\quad\quad\quad}
&
\;\;\;\;\;\;\;    \displaystyle   \lim_{d\to\infty}     K_*\big(   \cumpdxn  \big).
  }
\]

\par
We will devote the second half of the paper, section 5, 6 and 7, to prove the following result.
\par
{\bf Theorem 4.4.} {\it
Let $(X_n)_\nn$ be a sequence of finite metric spaces with uniform bounded geometry. If the coarse disjoint union $X=\bigsqcup_{n\in \mathbb{N}} X_n$ admits a fibred coarse embedding into
Hilbert space, then
\begin{small}
$$e_*:\; \lim_{d\to\infty}       K_*\big(   \culmpdxn  \big) \longrightarrow   \lim_{d\to\infty}     K_*\big(   \cumpdxn  \big)$$
\end{small}
is an isomorphism. Consequently, $\mu_{\max, \infty}$ is an isomorphism.
}
\hfill{$\Box$}
\par
We next prove Theorem 1.1 for the case $X=\bigsqcup_{n\in \mathbb{N}} X_n$ being a coarse disjoint union of a sequence of finite metric spaces, by using Theorem 4.4.
To do so, we need the following lemma (cf. \cite{GWY, OY, WiY2}).
\par
{\bf Lemma 4.5.} {\it
Let $X=\bigsqcup_{n\in \mathbb{N}} X_n$ be the coarse disjoint union of a sequence of finite metric spaces with uniform bounded geometry. For each $d\geq 0$, there is a short exact sequence
$$0\longrightarrow \cK\longrightarrow \cmpdx \longrightarrow \cumpdxn \longrightarrow 0$$
such that the inclusion $\cK\longrightarrow \cmpdx$ induces an injection on $K$-theory.
}
\par
{\bf Proof.} Let $Z_d\subset P_d(X)$ be a countable dense subset, and let $\zdn=Z_d\cap \pdxn$ for all $d\geq 0, \nn$ as being used in Definition 3.1 and Definition 4.1.
Note that $\cK\cong \cK(\ell^2(Z_d)\otimes H_0)$ is an ideal of $\cmpdx$.  There is a $*$-homomorphism
$$\Phi:\; \cpdx \longrightarrow \cupdxn$$
defined by $\Phi(T)=[(\Phi^{(0)}(T), \cdots, \Phi^{(n)}(T), \cdots)]$ for $T\in \cpdx,$
with
\[
\Phi^{(n)}(T) = \left\{
\begin{array}{ll}
0,                      & \mbox{ if } n<N_R; \\
T|_{\zdn\times \zdn},   & \mbox{ if } n\geq N_R,
\end{array}
\right.
\]
where $R=propagation(T)$ and $N_R\in \mathbb{N}$ is large enough such that
$$d\Big( X_n, \, \bigsqcup_{i=0}^{n-1}X_i \Big)>2R$$
for all $n\geq N_R$. The $*$-homomorphism $\Phi$ extends to $C^*_{\max}$-level
$$\Phi:\; \cmpdx \longrightarrow \cumpdxn$$
such that $\cK$ lives in the kernel of $\Phi$. The induced $*$-homomorphism on the quotient
$$\Phi:\; \cmpdx/\cK \longrightarrow \cumpdxn$$
has an inverse $\Psi$ defined as follows: for any
$$\ttt\in \cupdxn,$$
let $S=\bigoplus_\nn T^{(n)}$. Then $S\in \cpdx$ and $\Phi(S+\cK)=T$. Define
$$\Psi(T)=S+\cK\in \cmpdx/\cK.$$
Then $\Psi$ extends to a $*$-homomorphism on $C^*_{\max}$-level:
$$\Psi:\; \cumpdxn \longrightarrow \cmpdx/\cK$$
which is the inverse of $\Phi$. This gives the short exact sequence. The $K$-theory statement follows from \cite{GWY, OY}.
\hfill{$\Box$}

\par
{\bf Proposition 4.6.} {\it
Let $X=\bigsqcup_{n\in \mathbb{N}} X_n$  be the coarse disjoint union of a sequence of finite metric spaces with uniform bounded geometry. If
$X$ admits a fibred coarse embedding into Hilbert space, then the maximal coarse Baum-Connes conjecture holds for $X$. That is,
$$ \mu_{\max}: \lim_{d\to \infty} K_*(\pdx) \to  \lim_{d\to \infty} K_*(\cmpdx) \cong K_*(\cmx) $$
is an isomorphism.
}
\par
{\bf Proof.} For any $d\geq 0$, there exists $N_d\in \mathbb{N}$ large enough such that $d(X_n, X_m)>d$ provided $n, m\geq N_d$.  Let $X_{N_d}=\bigcup_{n=0}^{N_d-1} X_n$.
Then we have
$$K_*(\pdx)=K_*(P_d(X_{N_d}))\bigoplus \prod_{n=N_d}^\infty K_*(\pdxn).$$
By the definition of assembly maps and Lemma 4.5, we have the following commutative diagram:

\[
\xymatrix{ 0 \ar[d]& 0 \ar[d] \\
K_*(P_d(X_{N_d}))\oplus\bigoplus_{n=N_d}^\infty K_*(P_d(X_n)) \ar[d] \ar[r]& K_*(\mathcal{K}) \ar[d] \\
K_*(P_d(X)) \ar[d] \ar[r]^{\mu_{\max} \quad \quad} & K_*(C_{max}^*(\pdx)) \ar[d]  \\
\frac{\prod_{n=0}^\infty K_*(\pdxn)}{\oplus_{n=0}^\infty K_*(\pdxn)}
\ar[r]^{\mu_{\max, \infty}\quad\quad} \ar[d] & K_*\Big( \cumpdxn  \Big) \ar[d] \\ 0 & 0}
\]
\par
Passing to inductive limit as $d\to \infty$, the top horizontal arrow is an isomorphism for the following reason. An element in the sum, as a finite sequence, is supported on summands below some
fixed $m$ and, as $d\to \infty$, will eventually be absorbed into the first term on a single simplex. Thus, the assertion reduces to the fact that the assembly map is an isomorphism for a bounded metric
space. Now by Theorem 4.4 together with the five lemma, we complete the proof.
\hfill{$\Box$}
\par
Finally, we are able to prove Theorem 1.1 by using Proposition 4.6.
\par
{\bf Proof of Theorem 1.1.} Fix a point $x_0\in X$. For $n=0, 1, 2, \cdots$, let
$$X_n= \big\{ x\in X \; \big| \; n^3-n\leq d(x, x_0) \leq (n+1)^3+(n+1) \big\}$$
and denote
$$X^{(0)}=\bigcup_{n:\; even} X_n; \quad \quad X^{(1)}=\bigcup_{n: \; odd} X_n.$$
Then $X=X^{(0)}\cup X^{(1)}$, and each of $X^{(0)}$, $X^{(1)}$ and  $X^{(0)}\cap X^{(1)}$ is the coarse disjoint union of a sequence of finite subspaces of $X$, which admits a fibred coarse
embedding into Hilbert space as restricted from $X$.
\par
Moreover, the pair $(X^{(0)}, X^{(1)})$ is ``$\omega$-excisive'' (cf. \cite{HRY}) in the sense that, for any $R>0$ there exists $S>0$
such that
$$\Pen(X^{(0)}; R)\cap \Pen(X^{(1)}; R) \subseteq \Pen(X^{(0)}\cap X^{(1)}; S),$$
where $\Pen(Y; R)=\{x\in X|d(x, Y)\leq R\}$ is the $R$-neighborhood of a subspace. Indeed,
for any $R>0$, take an integer $n_R>R$ and let $S=(n_R+1)^3$. Suppose $x\in \Pen(X^{(0)}; R)\cap \Pen(X^{(1)}; R)$. If $d(x, x_0)\leq S$, then
$x\in \Pen(X^{(0)}\cap X^{(1)}; S)$ since $x_0\in X^{(0)}\cap X^{(1)}$. If $d(x, x_0)> S$, then there exist
$$\bar{x}\in X^{(0)}-B_X(x_0, n_R^3), \quad\quad \bar{y}\in X^{(1)}-B_X(x_0, n_R^3) $$
such that $d(x, \bar{x})\leq R$ and $d(x, \bar{y})\leq R$. Hence, $d(\bar{x}, \bar{y})\leq 2R$. We claim that either $\bar{x}\in X^{(0)}\cap X^{(1)}$, or $\bar{y}\in X^{(0)}\cap X^{(1)}$. Otherwise, we would simultaneously have that
$$\bar{x}\in \bigcup_{n:even,\, n\geq n_R} \big\{ x\in X\, \big|\, n^3+n<d(x, x_0)<(n+1)^3-(n+1) \big\}$$
and
$$\bar{y}\in \bigcup_{n:odd,\, n\geq n_R} \big\{ x\in X\, \big|\, n^3+n<d(x, x_0)<(n+1)^3-(n+1) \big\}.$$
It follows that $d(\bar{x}, \bar{y})\geq 2\, n_R>2R$, a contradiction. Consequently, either $\bar{x}$ or $\bar{y}$ is in $X^{(0)}\cap X^{(1)}$, so that
$x\in \Pen(X^{(0)}\cap X^{(1)}; R) \subset \Pen(X^{(0)}\cap X^{(1)}; S)$. Therefore, the pair $(X^{(0)}, X^{(1)})$ is ``$\omega$-excisive'' .
\par
For each $d\geq 0$, there exists $l_d>0$ such that, for any $x\in X\backslash B_X(x_0, l_d)$, the ball $B_X(x, d)$ is contained in either $X^{(0)}$ or $X^{(1)}$. Denote $K_d:=B_X(x_0, l_d)$. Then for Rips complexes, we have
$$P_d(X)=P_d(K_d\cup X^{(0)}) \bigcup P_d(K_d\cup X^{(1)}).$$
This is again an ``$\omega$-excisive'' decomposition,  into subspaces which are coarsely equivalent to $X^{(0)}$ and $X^{(1)}$ respectively.  As a result (cf. \cite{HRY}), we have the following
commutative diagram, in which the vertical arrows are assembly maps $\mu_{\max, \infty}$ connecting two Mayer-Vietoris exact sequences-----the top on $K$-homology, whereas the bottom on $K$-theory:

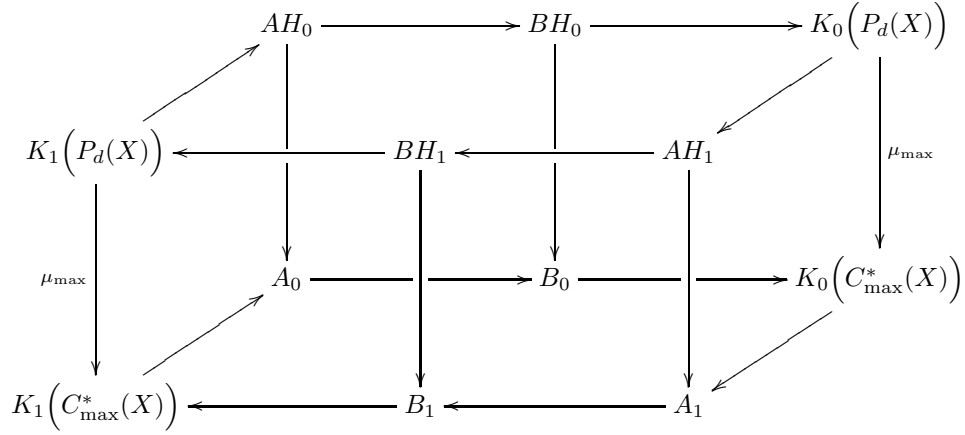
\[
\xymatrix
{
  & AH_0 \ar[rr] \ar'[d][dd]     &  &  BH_0\ar[rr]\ar'[d][dd]  &  &  K_0\Big(  \pdx \Big)  \ar[dd]^{\mu_{\max}} \ar[dl]        \\
    K_1\Big(  \pdx \Big) \ar[ur]\ar[dd]_{\mu_{\max}}    &  &  BH_1 \ar[ll]\ar[dd]     &  &  AH_1 \ar[ll]\ar[dd]          &    \\
  & A_0 \ar'[r][rr]              &  &  B_0 \ar'[r][rr]         &  &  K_0\Big(  \cmx  \Big) \ar[dl]                       \\
    K_1\Big(  \cmx  \Big) \ar[ur]                  &  &  B_1\ar[ll]              &  &  A_1 \ar[ll]                  &        }
\]
Where,
$$AH_0= K_0\Big(   P_d(K_d\cup X^{(0)}) \bigcap P_d(K_d\cup X^{(1)}) \Big) , \quad BH_0= K_0\Big(   P_d(K_d\cup X^{(0)})  \Big) \bigoplus   K_0\Big(   P_d(K_d\cup X^{(1)})  \Big) , $$
$$A_0=K_0\Big(  C^*_{\max} ( X^{(0)}\cap X^{1)})   \Big),  \quad\quad\quad\quad\quad\quad\quad\quad  B_0=K_0\Big(  C^*_{\max} ( X^{(0)}) \Big) \bigoplus K_0\Big(  C^*_{\max} ( X^{(1)})\Big) ,$$
and similarly for $K_1$. Passing to inductive limit as $d\to \infty$, clearly, Theorem 1.1 follows from (the proof of) Proposition 4.6 and the five lemma.
\hfill{$\Box$}

\section{The twisted algebras at infinity}
In the rest of this paper, we shall prove that the evaluation homomorphism
\begin{small}
$$e_*: \lim_{d\to \infty} K_*\Big( \culmpdxn \Big) \longrightarrow \lim_{d\to \infty} K_*\Big( \cumpdxn \Big) $$
\end{small}
is an isomorphism for a sequence of finite metric spaces $(X_n)_\nn$ with uniform bounded geometry such that the coarse disjoint union $X=\bigsqcup_\nn X_n$ admits a fibred coarse embedding
into Hilbert space. It follows from the discussion of the last section that this suffices to complete the proof of Theorem 1.1.
\par
In this section, we shall introduce the twisted Roe algebras at infinity  \\
$\cumpdxna$ and their localization counterpart \\
$\culmpdxna$ for $(X_n)_\nn$. Note that here the
coefficient algebras $\cA(V_n)$ are defined on {\em finite dimensional} affine subspaces $V_n$ of $H$. In section 6, we study various decompositions
for these twisted algebras to show that the evaluation map
\begin{small}
$$e_*: \lim_{d\to \infty} K_*\Big( \culmpdxna \Big) \longrightarrow \lim_{d\to \infty} K_*\Big( \cumpdxna \Big) $$
\end{small}
for the twisted algebras is an isomorphism. In section 7, we shall define the Bott maps $\beta$, $\beta_L$ and the Dirac maps $\alpha$, $\alpha_L$ to build the following
commutative diagram
\[
\xymatrix{
K_*\Big( \culmpdxn \Big) \ar[r]^{\quad e_*} \ar[d]^{(\beta_L)_*} & K_*\Big( \cumpdxn \Big) \ar[d]^{\beta_*} \\
K_*\Big( \culmpdxna \Big)  \ar@<1ex>[u]^{(\alpha_L)_*} \ar[r]^{\quad e_*}  & K_*\Big( \cumpdxna \Big )\;. \ar@<1ex>[u]^{\alpha_*}
}
\]
We prove a geometric analogue of the Bott periodicity in finite dimensions, i.e. $\alpha_*\circ \beta_*=identity$ and $(\alpha_L)_*\circ (\beta_L)_*=identity$. Passing to the inductive limit
as $d\to \infty$, a diagram chasing argument implies that the top evaluation map is an isomorphism,  as desired.
\par
\subsection*{{\it \S 5.1.  Preliminary}}
 Let $H$ be a separable Hilbert space. Denote by $V_a$, $V_b$ etc. the finite dimensional affine subspaces of $H$. Let $V_a^0$ be the linear subspace of $H$ consisting of differences of elements
of $V_a$. Let $\Cliff(V_a^0)$ be the complexified Clifford algebra of $V_a^0$ and $\cC(V_a)$ the graded $C^*$-algebra of continuous functions vanishing at infinity from $V_a$ into $\Cliff(V_a^0)$.
Let $\cS=C_0(\mathbb{R})$, graded according to even and odd functions. Define the graded tensor product:
$$\cA(V_a)=\cS\widehat{\otimes} \cC(V_a).$$
If $V_a\subseteq V_b$, then we have a decomposition $V_b=V_{ba}^0\oplus V_a$, where $V_{ba}^0$ is the orthogonal complement of $V_a^0$ in $V_b^0$. For each $v_b\in V_b$, we have a corresponding
decomposition $v_b=v_{ba}+v_a$, where $v_{ba}\in V_{ba}^0$ and $v_a\in V_a$. Every function $h$ on $V_a$ can be extended to a function $\widetilde{h}$ on $V_b$ by the formula $\widetilde{h}(v_{ba}+v_a)=h(v_a)$.
\par
{\bf Definition 5.1.} For affine subspaces $V_a\subseteq V_b$, denote by $C_{ba}: V_b\to \Cliff(V_{ba}^0)$ the function $v_b\mapsto v_{ba}\in \Cliff(V_{ba}^0)$, where $v_{ba}$ is considered as an element
of $\Cliff(V_{ba}^0)$ via the inclusion $V_{ba}^0\subseteq \Cliff(V_b^0)$. Let $X$ be the unbounded multiplier of $\cS$ given by the function $t\mapsto t$. Define a $*$-homomorphism
$\beta_{V_b, V_a}:\cA(V_a)\to \cA(V_b)$, or simply denoted by $\beta_{ba}$, by the formula
$$\beta_{ba}(g\widehat\otimes h)=g(X\widehat\otimes 1+ 1\widehat\otimes C_{ba})\, (1\widehat\otimes \wt h) $$
for all $g\in \cS, h\in \cC(V_a)$, where $g(X\widehat\otimes 1+1\widehat\otimes C_{ba})$ is defined by functional calculus.
\par
{\bf Definition 5.2.} If  $V_a\subseteq V_b$, for any subset $O\subset \mathbb{R}_+\times V_a$, define
$$\overline{O}^{\beta_{ba}}=\Big\{ \; \big(t, \;\; v_{ba}+v_a \big) \in \mathbb{R}_+\times V_b \; \Big| \;\; \Big( \sqrt{t^2+\|v_{ba}\|^2}\; ,\;  v_a \Big)\in O \; \Big\}.$$
\par
For any finite dimensional affine subspace $V_a$ of $H$, the algebra $C_0(\mathbb{R}_+\times V_a)$ is included in $\cA(V_a)$ as its center. For any function $a\in \cA(V_a)$, the support of $a$,
denoted by $\supp(a)$, is the complement of all points $(t, v)\in \mathbb{R}_+\times V_a$ such that there exists $g\in C_0(\mathbb{R}_+\times V_a)$ such that $g(t, v)\not=0$ but $g\cdot a=0$. Note
that if $V_a\subseteq V_b$ and $a\in \cA(V_a)$, then
$$\supp\big( \beta_{ba}(a) \big) = \overline{\supp(a)}^{\beta_{ba}}.$$
\par
{\bf Definition 5.3.} Let $W_a, W_b, V_a, V_b, V_c$ be finite dimensional affine subspaces of $H$ with  $W_a\subseteq W_b\subset V_c$ and $V_a\subseteq V_b\subset V_c$.
Let $t$ be an affine isometry from $W_b$ onto $V_b$ mapping $W_a$ onto $V_a$. Then $t$ canonically
induces a $*$-isomorphism from $\cA(W_a)$ onto $\cA(V_a)$. Note that $\cA(V_a)$ is included in $\cA(V_c)$ via $\beta_{ca}$. We denote by $t_*$ the composition
$$t_*: \;\; \cA(W_a)\stackrel{\cong}{\longrightarrow} \cA(V_a) \stackrel{\beta_{ca}}{\longrightarrow} \cA(V_c).$$
Then we have the following commutative diagram which is useful in the definition of product structure for the twisted Roe algebras.
\[
\xymatrix{
t_*: & \cA(W_a) \ar[r]^{\cong} \ar[d]^{\beta_{ba}} & \cA(V_a) \ar[r]^{\beta_{ca}} \ar[d]^{\beta_{ba}}  & \cA(V_c) \ar[d]^{=}   \\
& \cA(W_b) \ar[r]^{\cong}                     & \cA(V_b) \ar[r]^{\beta_{ca}} & \cA(V_c)
}
\]
\hfill{$\Box$}
\par
\subsection*{\it \S 5.2. The twisted Roe algebras at infinity}
Let $(X_n)_\nn$ be a sequence of finite metric spaces with uniform bounded geometry such that the coarse disjoint union $X=\bigsqcup_\nn X_n$ admits a fibred coarse embedding into Hilbert space.
To fix notations, we note that the notion of fibred coarse embedding into Hilbert space for the coarse disjoint union $X$ is equivalent to the following notion for the sequence
$(X_n)_\nn$:
\par
{\bf Definition 5.4.} A sequence of finite metric spaces $(X_n)\nn$ with uniform bounded geometry is said to admit a fibred coarse embedding into Hilbert space if there exist
\begin{itemize}
\item a field of Hilbert spaces $(H_x)_{x\in X_n, \nn}$;
\item a section $s: X_n\to\displaystyle\bigsqcup_{x\in X_n} H_x$ for all $\nn$;
\item two non-decreasing functions $\rho_1$ and $\rho_2$ from $[0, \infty)$ to $(-\infty, \infty)$ with $\lim_{r\to \infty}\rho_i(r)=\infty$ ($i=1, 2$);
\item a non-decreasing sequence of numbers $0\leq \l_0\leq \l_1\leq \cdots \leq \l_n\leq \cdots$ with $\lim_{n\to \infty} \l_n=\infty$
\end{itemize}
such that for each $x\in X_n, \nn$ there exists a trivialization
$$t_x: (H_z)_{z\in B(x, l_n)}\longrightarrow B(x, l_n)\times H$$
such that the restriction  of $t_x$ to the fiber $H_z$ ($z\in B(x, l_n)$) is an affine isometry $t_x(z): H_z\to H$, satisfying
\begin{itemize}
\item[(1)] $\rho_1(d(z, z'))\leq \|t_x(z)(s(z))-t_x(z')(s(z'))\|\leq \rho_2(d(z, z'))$ for any $z, z'\in B(x, l_n)$, $x\in X_n, \nn$;
\item[(2)] for any $x, y\in X_n$ with $B(x, l_n)\cap B(y, l_n)\not=\emptyset$, there exists an affine isometry $t_{xy}: H\to H$ such that
		$t_{x}(z)\circ t_{y}^{-1}(z)=t_{xy}$ for all $z\in B(x, l_n)\cap B(y, l_n)$.
\end{itemize}
\hfill{$\Box$}
\par
For each $d\geq 0$ and $\nn$, let $\pdxn$ be the Rips complex of $X_n$ at scale $d$ endowed with the spherical metric. For each $x\in X_n$, denote by $\Star(x)$ the open star of $x$ in
the second barycentric subdivision of $\pdxn$. Take a countable dense subset $\zdn\subset \pdxn$ for each $d\geq 0$ in such a way that
$$ (1) \quad \zdn\subset \displaystyle\bigsqcup_{x\in X_n} \Star(x);  \quad\quad\quad (2) \quad \zdn\subseteq Z_{d', n}  \mbox{ when }  d<d'. $$
\par
For any $x\in \zdn$, there exists a unique point $\bar{x}\in X_n$ such that $x\in \Star(\bar{x})$. We define
$$H_x=H_{\bar{x}}, \quad\quad s(x)=s(\bar{x})$$
for all $x\in \zdn\cap \Star(\bar x)$ and let
$$t_x(z)=t_{\bar{x}}(\bar{z}):\; H_z\to H$$
for all $x, z\in \zdn, \; \nn$ with $\bar z\in B(\bar x, l_n)$.
\par
For each $\nn$, define $V_n$ to be the {\em finite dimensional} affine subspace of $H$ spanned by $t_x(z)(s(z))$ for all $z\in B(x, l_n), \; x\in X_n$:
$$V_n :=\mbox{ affine-span } \Big\{ \; t_x(z)(s(z))\; \Big| \; z\in B(x, l_n), \; x\in X_n \Big\} \,.$$
For each $x\in \zdn$, $k\geq 0$, define
$$W_k(x)=\mbox{ affine-span } \Big\{ \; t_x(z)(s(z))\; \Big| \; z\in B_{\pdxn} (x, k)  \; \Big\} \subseteq V_n \,.$$
Note that for each $k\geq 0$, there exists $N\in \mathbb{N}$ such that  $W_k(x)$ is well-defined for $x\in \zdn$ with $n\geq N$, and is an affine subspace of $V_n$. By Definition 5.1,
the inclusion $W_k(x)\subseteq V_n$ induces the inclusion
$$\beta_{V_n, W_k(x)}:\; \cA(W_k(x))\to \cA(V_n).$$
\par
{\bf Definition 5.5.} For each $d\geq 0$, define $\cupdxna$ to be the set of all equivalence classes $\ttt$ of sequences $\big( T^{(0)}, \cdots, T^{(n)}, \cdots \big)$ described as follows:
\begin{itemize}
\item[(1)] $T^{(n)}$ is a bounded function from $\zdn \times \zdn$ to $\cA(V_n)\widehat\otimes \cK$ for all $\nn$;
\item[(2)] for any bounded subset $B\subset \pdxn$, the set
			$$\{ (x, y)\in B\times B \cap \zdn\times \zdn \; | \; T^{(n)} (x, y)\not= 0 \}$$
		is finite;
\item[(3)] there exists $L>0$ such that
			$$\#\{y\in \zdn \; | \; T^{(n)}(x, y)\not= 0 \}<L, \quad \#\{y\in \zdn \; | \; T^{(n)}(y, x)\not= 0 \}<L $$
			for all $x\in \zdn,\; \nn$;
\item[(4)] there exists $R>0$ such that $T^{(n)}(x, y)=0$ whenever $d(x, y)>R$  for $x, y\in \zdn,\; \nn$; (the least such $R$ is call the propagation of the sequence
			 $\big( T^{(0)}, \cdots, T^{(n)}, \cdots \big)$ .)
\item[(5)] there exists $r>0$ such that
				\[
					\begin{array}{rcl}
						\supp(T^{(n)}(x, y) )  & \subseteq & B_{\mathbb{R}_+\times V_n} \big( t_x(x)(s(x)), \; r\big)     \\
										  & :=   & \big\{ (\tau, v)\in \mathbb{R}_+\times V_n \; \big| \; \tau^2+\|v-t_x(x)(s(x))\|^2<r^2 \big\}
					\end{array}
				\]
		for all $x, y\in \zdn, \; \nn$.
\item[(6)] there exist $k\geq 0$ and $K>0$ depending only on the sequence $\big( T^{(0)}, \cdots, T^{(n)}, \cdots \big)$ (not on $n$) such that, for each $x, y\in \zdn$, there exists
			$$T_1^{(n)}(x, y)\in \cA(W_k(x))\widehat\otimes \cK\cong \cS\widehat\otimes \cC(W_k(x))\widehat\otimes \cK$$
			of the form $\sum_{i=1}^K g_i\widehat\otimes h_i\widehat\otimes k_i$ where $g_i\in \cS$, $h_i\in \cC(W_k(x))$,
			$k_i\in \cK$ for $1 \leq i \leq K$  such that
				$$ T^{(n)}(x, y) = \Big( \beta_{V_n, W_k(x)} \widehat\otimes 1 \Big) \big( T^{(n)}_1 (x, y) \big) \,.$$
\item[(7)] there exists $c>0$ such that if $T_1^{(n)}(x, y)\in \cA(W_k(x))\widehat\otimes \cK$  as above, and $w\in \mathbb{R}_+\times W_k(x)$ is of norm one, then
		the derivative of $T^{(n)}_1(x, y)$ in the direction $w$, $\nabla_w\big( T^{(n)}_1(x, y) \big)$, exists in $\cA(W_k(x))\widehat\otimes \cK$ and is of norm at most $c$.
\end{itemize}
The equivalence relation $\sim$ on these sequences is defined by
$$(T^{(0)}, \cdots, T^{(n)}, \cdots) \sim (S^{(0)}, \cdots, S^{(n)}, \cdots)$$
if and only if
$$\lim_{n\to \infty} \sup_{x, y\in \zdn} \big\| T^{(n)}(x, y) - S^{(n)}(x, y) \big\|_{\cA(V_n)\widehat\otimes \cK} =0.$$
The product structure for $\cupdxna$ is defined as follows. For any two elements $\ttt$ and $S=\big[ (S^{(0)}, \cdots, S^{(n)}, \cdots) \big]$ in $\cupdxna$, their product is defined to be
$$\big[ \big( (TS)^{(0)}, \cdots, (TS)^{(n)}, \cdots \big) \big] \,,$$
where there exists a sufficiently large $N\in \mathbb{N}$ depending on the propagation of the two representative sequences, such that
$ (TS)^{(n)}=0$ for all $n=0, 1, 2, \cdots, N-1$ and
$$(TS)^{(n)}(x, y)=\sum_{z\in \zdn} \Big( T^{(n)}(x, z) \Big) \cdot \Big( \big(t_{xz}\big)_* \big( S^{(n)}(z, y) \big) \Big)$$
for all $x, y\in \zdn, \; n\geq N$.
\par
{\bf Remark 5.6.} Some explanations to the above formula are given here. Note that there exists $k\geq 0$ such that
				$ S^{(n)}(z, y) = \Big( \beta_{V_n, W_k(z)} \widehat\otimes 1 \Big) \big( S^{(n)}_1 (z, y) \big)$ for $z, y\in \zdn, \;\nn$, where
				$S_1^{(n)}(z, y)\in \cA(W_k(z))\widehat\otimes \cK$.  For the propagation $R$ of the representative sequence $\big( T^{(0)}, \cdots, T^{(n)}, \cdots \big)$ and the above $k$, there
exists $N\in \mathbb{N}$ large enough such that for any $n\geq N$, there exists $\widetilde{l_n}>0$ (since $\pdxn$ is coarsely equivalent to $X_n$)
such that the trivialization
$$t_x: \;\; \Big\{ H_z \Big\}_{z\in B_{\pdxn}(x, \widetilde{l_n})} \longrightarrow  B_{\pdxn}(x, \widetilde{l_n})\times H$$
makes sense on $B_{\pdxn}(x, \widetilde{l_n})$, and $R+k<\widetilde{l_n}$. If $d(x, z)\leq R$ in $\pdxn$, then the affine isometry $t_{xz}=t_x(w)\circ t_z^{-1}(w):\; H\to H$ for all
$$w\in B_{\pdxn}(x, \widetilde{l_n})\cap B_{\pdxn}(z, \widetilde{l_n}).$$
Note that $t_{xz}$ maps the affine subspace
$$W_k(z)=\mbox{ affine-span } \Big\{ t_z(w)(s(w)) \; \Big| \; w\in B_{\pdxn}(z, k) \cap \zdn \Big\}$$
onto an affine subspace of $W_{R+k}(x)$ since $B_{\pdxn}(z, k)\subseteq B_{\pdxn}(x, R+k)$. By Definition 5.3, the composition
$$t_{xz}: W_k(z)\to t_{xz}\big( W_k(z) \big) \subseteq W_{R+k}(x) \subseteq V_n$$
induces the $*$-homomorphism
$$\big(t_{xz}\big)_*:\; \cA( W_k(z) ) \to \cA(V_n).$$
We define
$$\big(t_{xz}\big)_*\Big( S^{(n)}(z, y) \Big):= \big(t_{xz}\big)_*\Big( S_1^{(n)}(z, y) \Big)$$
in the above product formula. Note that for $\nn$ large enough, this definition of $\big(t_{xz}\big)_*\Big( S^{(n)}(z, y) \Big)$ does not depend on the choice of $k$
(see Definition 5.3, and \cite{HG, HKT}).
\hfill{$\Box$}
\par
The $*$-structure for $\cupdxna$ is defined by the formula
$$ \big[ \big( T^{(0)}, \cdots, T^{(n)}, \cdots \big) \big]^*=   \big[ \big( (T^*)^{(0)}, \cdots, (T^*)^{(n)}, \cdots \big) \big] \,,$$
where
$$(T^*)^{(n)}(x, y)=\big( t_{xy} \big)_* \Big( \big( T^{(n)}(y, x)\big)^*\Big)$$
for all but finitely many $n$, and $0$ otherwise.
\par
Now, $\cupdxna$ is made into a $*$-algebra by using the additional usual matrix operations. Define $\cumpdxna$ to be the completion of $\cupdxna$ with respect to the norm
\begin{footnotesize}
$$\big\|T\big\|_{\max} :=\sup\Big\{ \big\|\phi(T)\big\|_{\cB(H_\phi)} \; \Big | \; \phi:\cupdxna \to \cB(H_\phi), \mbox{ a $*$-representation} \Big\}. $$
\end{footnotesize}
\hfill{$\Box$}
\par
{\bf Definition 5.7.}  For each $d\geq 0$, define $\culpdxna$  to be the $*$-algebra of all bounded and uniformly norm-continuous functions
$$f: \; [0,\infty)\longrightarrow \cupdxna$$
such that $f(t)$ is of the form $f(t)=[(f^{(0)}(t), \cdots, f^{(n)}(t), \cdots )]$ for all $t\in [0, \infty)$, where the family of functions
$(f^{(n)}(t))_{\nn, t\geq 0}$ satisfy the conditions in Definition 5.5 with uniform constants, and there exists a bounded function $R:\; [0, \infty)\to [0,\infty)$ with
$\lim_{t\to \infty} R(t)=0$ such that
$$ \big( f^{(n)}(t)\big) (x, y) =0  \quad \mbox{whenever} \quad d(x, y)>R(t)$$
for all $x, y\in \zdn$, $\nn$, $t\in [0, \infty)$.
\par
Define $\culmpdxna$ to be the completion of
$$\culpdxna$$
with respect to the norm
$$\big\| f \big\|_{\max} :=\sup_{t\in [0, \infty)} \big\|f(t)\big\|_{\max}.$$
\hfill{$\Box$}
\par
The evaluation homomorphism
$$e: \; \culmpdxna \longrightarrow \cumpdxna $$
defined by $e(f)=f(0)$ induces the evaluation homomorphism on $K$-theory:
\begin{small}
$$ \lim_{d\to\infty}       K_*\big(   \culmpdxna  \big) \stackrel{e_*}{\longrightarrow}   \lim_{d\to\infty}     K_*\big(   \cumpdxna  \big).$$
\end{small}

\section{Decompositions of twisted algebras}
In this section we shall prove the following result.
\par
{\bf Theorem 6.1.} {\it
Let $(X_n)_\nn$ be a sequence of finite metric spaces with uniform bounded geometry which admit a fibred coarse embedding into Hilbert space. Then
the evaluation homomorphism
\begin{small}
$$ \lim_{d\to\infty}       K_*\big(   \culmpdxna  \big) \stackrel{e_*}{\longrightarrow}   \lim_{d\to\infty}     K_*\big(   \cumpdxna  \big).$$
\end{small}
is an isomorphism.
}
\par
The proof proceeds by decomposing the twisted algebras into various smaller subalgebras and applying a Mayer-Vietoris argument.
\subsection*{\it \S 6.1. Coherent system of open subsets of $\mathbb{R}_+\times V_n$ ($\nn$)}
To begin with, we shall discuss ideals of the twisted algebras supported on certain open subsets.
\par
{\bf Definition 6.2.} A collection $O=(O_{n, x})_{x\in X_n, \,\nn}$ of open subsets of $\mathbb{R}_+\times V_n, \; \nn,$ is said to be a {\em coherent system} if for all but finitely many $\nn$,
the following conditions hold:
\begin{itemize}
\item[(1)] for any subset $C\subseteq B(x, l_n)\cap B(y, l_n)$ with $x, y\in X_n$, we have
		$$ O_{n, x}\cap \big(\mathbb{R}_+\times W_C(x) \big) = t_{xy} \Big( O_{n, y} \cap \big(\mathbb{R}_+\times W_C(y) \big) \Big),$$
		where
		$$W_C(x)=\mbox{ affine-span } \Big\{ t_x(z)(s(z)) \; \Big| \; z\in C \Big\}=t_{xy}\big( W_C(y) \big),$$
		$$W_C(y)=\mbox{ affine-span } \Big\{ t_y(z)(s(z)) \; \Big| \; z\in C \Big\}=t_{yx}\big( W_C(x) \big),$$
		and $t_{xy}=t_x(z)\circ t_y^{-1}(z):\; H\to H$ for all $z\in B(x, l_n)\cap B(y, l_n)$.
\item[(2)] for any $C\subseteq B(x, l_n)$, where $x\in X_n, \;\nn$, and any affine subspace $W$ with $W_C(x)\subseteq W\subseteq V_n$, we have
		$$ \overline{O_{n, x}\cap \big(\mathbb{R}_+\times W_C(x) \big)}^{\beta_{W, W_C(x)}} \subseteq  O_{n, x}\cap \big(\mathbb{R}_+\times W \big),$$
		where recall (Definition 5.3) that for $V_a\subseteq V_b$ and $O\subseteq \mathbb{R}_+\times V_a$,
		$$ \overline{O}^{\beta_{ba}} = \big\{  (t, v_{ba}+v_a)\in \mathbb{R}_+\times V_b \; \big| \; \big( \sqrt{t^2+\|v_{ba}\|^2} \; , \;v_a \big)\in O \big\}.$$
\end{itemize}
 \hfill{$\Box$}
\par
{\bf Definition 6.3.} Let $O=(O_{n, x})_{x\in X_n, \,\nn}$ be a coherent system of open subsets of $\mathbb{R}_+\times V_n, \; \nn$. For each $d\geq 0$, define
$$\cupdxna_O$$
to be the $*$-subalgebra of $\cupdxna$ generated by the equivalence classes of those sequences $\big[ \big( T^{(0)}, \cdots, T^{(n)}, \cdots \big) \big]$ such that
$$\supp (T^{(n)}(x, y)) \subseteq O_{n, \bar{x}}$$
for all $x, y\in \zdn$ with $x\in \Star(\bar x)$ for all $n\geq N$ for some $n\in \mathbb{N}$ large enough depending on the sequence  $\big[ \big( T^{(0)}, \cdots, T^{(n)}, \cdots \big) \big]$.
\par
Define $\cumpdxna_O$ to be the norm closure of  \\
$\cupdxna_O$ in $\cumpdxna$.
 \hfill{$\Box$}
\par
{\bf Lemma 6.4.} {\it
$\cumpdxna_O$ is a two sided ideal of    \\
$\cumpdxna$.
}
\par
{\bf Proof.} Suppose
$$\big[ \big( T^{(0)}, \cdots, T^{(n)}, \cdots \big) \big]\in \cupdxna  \,,$$
$$\big[ \big( S^{(0)}, \cdots, S^{(n)}, \cdots \big) \big]\in \cupdxna_O  \,,$$
so that (1) there exists $R>0$ such that the propagation of $T^{(n)}$ and $S^{(n)}$ is at most $R$ for all $\nn$; (2) $\supp(S^{(n)}(z, y))\subseteq O_{n,z}$ for all $z, y\in \zdn, \,\nn$;
(3) there exists $k\geq 0$ such that $S^{(n)}(z, y)=\big( \beta_{V_n, W_k(z)}\widehat\otimes 1\big) \big( S_1^{(n)}(z, y) \big)$ for some $S_1^{(n)}(z, y)\in \cA(W_k(z))\widehat\otimes \cK$
for all $z, y \in \zdn \,\nn$. It follows that
$$\supp \big(S_1^{(n)}(z, y) \big) \subseteq O_{n, z}\cap \big( \mathbb{R}_+\times W_k(z) \big), $$
and
\[
\begin{array}{rcl}
\supp\Big( \big( t_{xz} \big)_* \big( S_1^{(n)}(z, y) \big) \Big) & \subseteq & t_{xz} \big(  O_{n, z}\cap (\mathbb{R}_+\times W_k(z)) \big)     \\
  & = & O_{n, x}\cap \big( \mathbb{R}_+\times W_{B(z, k)}(x) \big)  \,,
\end{array}
\]
where $B(z, k):=B_{\pdxn}(z, k)$ and
$$W_k(z)=\mbox{ affine-span }\Big\{ t_z(w)(s(w)) \; \Big| \; w\in B_{\pdxn}(z, k) \cap \zdn \Big\},$$
$$W_{B(z, k)}(x) = t_{xz}(W_k(z)).$$
It follows from the definition of a coherent system that
\[
\begin{array}{rcl}
\supp\Big( \big( t_{xz} \big)_* \big( S^{(n)}(z, y) \big) \Big) & \subseteq & \overline{ \Big(  O_{n, x}\cap (\mathbb{R}_+\times W_{B(z, k)}(x)) \Big) } ^{\beta_{V_n, W_{B(z, k)}(x)}}    \\
  & \subseteq  & O_{n, x} \,.
\end{array}
\]
Hence,
$$\supp\Big( T^{(n)}(x, z) \cdot \big( t_{xz} \big)_* \big( S^{(n)}(z, y) \big) \Big) \subseteq O_{n, x}, $$
for all but finitely many $\nn$. That is,
$$\big[ \big( (TS)^{(0)}, \cdots, (TS)^{(n)}, \cdots \big) \big]\in \cupdxna_O.$$
The proof is complete.
\hfill{$\Box$}
\par
{\bf Definition 6.5.} Let $O=(O_{n, x})_{x\in X_n, \,\nn}$ be a coherent system of open subsets of $\mathbb{R}_+\times V_n, \; \nn$. For each $d\geq 0$, define
$$\culpdxna_O$$
to be the $*$-subalgebra of $\culpdxna$ consisting of all functions
$$f:\; [0, \infty) \longrightarrow \cupdxna_O.$$
Define $\culmpdxna_O$ to be the completion of    \\
$\culpdxna_O$ in $\culmpdxna$.
\hfill{$\Box$}
\par
Note that $\culmpdxna_O$ is an ideal of    \\
$\culmpdxna$. We also have an  evaluation homomorphism
$$e: \; \culmpdxna_O\rightarrow \cumpdxna_O$$
given by $e(f)=f(0)$.
\par

\subsection*{\it \S  6.2. Strong Lipschitz homotopy invariance}
In this subsection, we shall investigate ideals of the twisted algebras supported on those coherent systems which are separated by subsets of $X_n$, $\nn$. These ideals can be decomposed into certain algebras
defined on uniformly bounded subsets of $P_d(X_n)$, $\nn$, whose $K$-theory turns out to be strongly Lipschitz homotopy invariant. As a result, the evaluation homomorphism for these ideals is
an isomorphism. This  constitutes a major step towards the proof of Theorem 6.1.
\par
Let $\Gamma_n$ be a subset of $X_n$ for each $\nn$ and denote $\Gamma=(\Gamma_n)_\nn$. Let $r>0$.
\par
{\bf Definition 6.6.} A coherent system $O=(O_{n, x})_{x\in X_n, \,\nn}$ of open subsets of  $\mathbb{R}_+\times V_n, \; \nn,$ is said to be {\em $(\Gamma, r)$-separate} if
there exist open subsets $\big( O_{n, x, \gamma} \big)_{\gamma\in \Gamma_n\cap B(x, l_n)}$ of $\mathbb{R}_+\times V_n$ for all $x\in X_n,\, \nn$, such that
\begin{itemize}
\item $O_{n, x}=\displaystyle\bigcup_{\gamma\in \Gamma_n\cap B(x, l_n)} O_{n, x, \gamma} $ ;
\item $ O_{n, x, \gamma} \cap O_{n, x, \gamma'} =\emptyset $ for distinct $\gamma, \gamma'\in \Gamma_n\cap B(x, l_n)$ ;
\item $ O_{n, x, \gamma} \subseteq B_{\mathbb{R}_+\times V_n}\Big( t_x(\gamma)(s(\gamma)) \, ,  \, r \Big)$ for each $\gamma\in \Gamma_n\cap B(x, l_n)$.
\end{itemize}
\hfill{$\Box$}
\par
{\bf Theorem 6.7.} {\it
If a coherent system $O=(O_{n, x})_{x\in X_n, \,\nn}$ is  $(\Gamma, r)$-separate, then the evaluation homomorphism on $K$-theory
\[
\begin{array}{ll}
e_*: &   \displaystyle \lim_{d\to \infty} K_*\Big(  \culmpdxna_{O} \Big)             \\
&   \longrightarrow  \displaystyle \lim_{d\to \infty} K_*\Big(  \cumpdxna_{O} \Big)
\end{array}
\]
is an isomorphism.
}
\par
We need some preparations before we can prove Theorem 6.8. Suppose $O=(O_{n, x})_{x\in X_n, \,\nn}$ is  $(\Gamma, r)$-separate for some  $\Gamma=(\Gamma_n)_\nn$ and $r>0$ as above.
For $d\geq 0$, let $(Y_\gamma)_{\gamma\in \Gamma_n}$ be a collection of closed subsets of $\pdxn$ for each $\nn$ such that
\begin{itemize}
\item $(Y_\gamma)_{\gamma\in \Gamma_n, \,\nn}$ is uniformly bounded, i.e. there exists $K>0$ such that $diameter(Y_\gamma)\leq K$ for all $\gamma\in \Gamma_n, \nn$.
\item $\gamma\in Y_\gamma$ for all $\gamma\in \Gamma_n, \nn$.
\end{itemize}
In particular, we are mainly concerned with the following cases:
\begin{itemize}
\item[(1)] $Y_\gamma=B(\gamma, S):=\big\{ x\in \pdxn \;\big| \; d(x, \gamma)\leq S \big\}$ for some common $S>0$ for all $\gamma\in \Gamma_n, ,\nn$.
\item[(2)] $Y_\gamma=\Delta_\gamma(S)$, the simplex in $\pdxn$ with vertices $\{x\in X_n | d(x, \gamma)\leq S \}$ for some common $S>0$ for all $\gamma\in \Gamma_n, \,\nn$.
\item[(3)] $Y_\gamma=\{\gamma\}$ for all $\gamma\in \Gamma_n, \,\nn$.
\end{itemize}
\par
{\bf Definition 6.8.} For an open subset $O\subseteq \mathbb{R}_+\times V_n$, we denote by $\cA(O)$ the $C^*$-subalgebra of $\cA(V_n)$ generated by the functions whose supports are
contained in $O$. Note that $\cA(O)$ is an ideal of $\cA(V_n)$.
\par
{\bf Definition 6.9.}  Define $A_\infty[ \, ( Y_\gamma: \; \gamma\in \Gamma_n)_\nn]$ to be the $*$-subalgebra of
\begin{equation}\label{qqq}
\frac
{\prod_\nn \Big( \bigoplus_{\gamma\in \Gamma_n}  C^*_{\max}(Y_\gamma)\widehat\otimes \cA(O_{n, \gamma, \gamma}) \Big) }
{\bigoplus_\nn \Big( \bigoplus_{\gamma\in \Gamma_n}  C^*_{\max}(Y_\gamma)\widehat\otimes \cA(O_{n, \gamma, \gamma}) \Big) }
\end{equation}
consisting of elements of the form $\big[ \big( T^{(0)}, \cdots, T^{(n)},  \cdots \big) \big]$ where
$$T^{(n)}=\bigoplus_{\gamma_\in \Gamma_n} T^{(n)}_\gamma$$
with
$$T^{(n)}_\gamma\in \mathbb{C}[Y_\gamma]\widehat\otimes \cA(O_{n, \gamma, \gamma}) \subset C^*_{\max}(Y_\gamma)\widehat\otimes \cA(O_{n, \gamma, \gamma}),$$
and when viewed as functions
$$T^{(n)}_\gamma:\; (\zdn\times \zdn)\cap (Y_\gamma\times Y_\gamma) \rightarrow \cA(O_{n, \gamma, \gamma}),$$
the family  $(T^{(n)}_\gamma)_{\gamma\in \Gamma_n, \nn}$ satisfy the conditions in Definition 5.5 with uniform constants.
\par
Define $A^*_\infty ( \, ( Y_\gamma: \; \gamma\in \Gamma_n)_\nn )$ to be the completion of
$$A_\infty[ \, ( Y_\gamma: \; \gamma\in \Gamma_n)_\nn]$$
inside the $C^*$-algebra (1).
\hfill{$\Box$}
\par
{\bf Definition 6.10.} Define $A_{L, \infty}[ \, ( Y_\gamma: \; \gamma\in \Gamma_n)_\nn]$ to be the $*$-algebra of all bounded and uniformly norm-continuous functions
$$f: \; [0,\infty)\longrightarrow A_\infty[ \, ( Y_\gamma: \; \gamma\in \Gamma_n)_\nn]$$
where $f(t)$ is of the form $f(t)=[(f^{(0)}(t), \cdots, f^{(n)}(t), \cdots )]$ for all $t\in [0, \infty)$ and
$$f^{(n)}(t)=\bigoplus_{\gamma_\in \Gamma_n} f^{(n)}_\gamma (t)$$
such that the family of functions
$(f_\gamma^{(n)}(t))_{\gamma\in \Gamma_n, \nn, t\geq 0}$ satisfy the conditions in Definition 5.5 with uniform constants, and there exists a bounded function $R:\; [0, \infty)\to [0,\infty)$ with
$\lim_{t\to \infty} R(t)=0$ such that
$$ \big( f_\gamma^{(n)}(t)\big) (x, y) =0  \quad \mbox{whenever} \quad d(x, y)>R(t)$$
for all $x, y\in \zdn\cap Y_\gamma$, $\gamma\in \Gamma_n$, $\nn$, $t\in [0, \infty)$.
\par
Define $A^*_{L, \infty} \big( \, ( Y_\gamma: \; \gamma\in \Gamma_n)_\nn \big)$ to be the completion of
$$A_{L, \infty}[ \, ( Y_\gamma: \; \gamma\in \Gamma_n)_\nn]$$
with respect to the norm
$$\big\| f \big\|_{\max} :=\sup_{t\in [0, \infty)} \big\|f(t)\big\|_{\max}.$$
\hfill{$\Box$}
\par
{\bf Proposition 6.11.} {\it
Suppose a coherent system $O=(O_{n, x})_{x\in X_n, \,\nn}$ is  $(\Gamma, r)$-separate for some  $\Gamma=(\Gamma_n)_\nn$ and $r>0$ as above. Then
\begin{itemize}
\item  $\displaystyle  \cumpdxna_O \cong \lim_{S\to \infty} A^*_{\infty}\big( \,\big( B(\gamma, S);  \; \gamma\in \Gamma_n  \big)_\nn \big)  \, ,$
\item  $\displaystyle  \culmpdxna_O \cong \lim_{S\to \infty} A^*_{L, \infty}\big( \,\big( B(\gamma, S);  \; \gamma\in \Gamma_n  \big)_\nn \big)  \, ,$
\item  $\displaystyle \lim_{d\to \infty}\cumpdxna_O \cong \lim_{S\to \infty} A^*_{\infty}\big( \,\big( \Delta_\gamma (S);  \; \gamma\in \Gamma_n  \big)_\nn \big)  \, ,$
\item  $\displaystyle  \lim_{d\to \infty}\culmpdxna_O \cong \lim_{S\to \infty} A^*_{L, \infty}\big( \,\big( \Delta_\gamma (S);  \; \gamma\in \Gamma_n  \big)_\nn \big)  \,. $
\end{itemize}
}
\par
{\bf Proof.} We shall establish an isomorphism for the first item. Take arbitrarily an element
$$\ttt\in \cupdxna_O,$$
where the functions $T^{(n)}:\; \zdn\times \zdn\rightarrow \cA(V_n)\widehat\otimes \cK$ have supports
$$\supp\big( T^{(n)}(x, y) \big) \subseteq O_{n, x} \subseteq \bigsqcup_{\gamma\in \Gamma_n\cap B_{X_n}(x, l_n)} O_{n, x, \gamma}$$
for all $x, y\in \zdn, \nn, $ since the coherent system $O$ is $(\Gamma, r)$-separate. Then we have a direct sum decomposition
$$T^{(n)}(x, y)=\bigoplus_{\gamma\in \Gamma_n\cap B_{X_n}(x, l_n)}  T_\gamma^{(n)}(x, y) \,,$$
where
$$T_\gamma^{(n)}(x, y) = T^{(n)}(x, y)\Big|_{O_{n, x, \gamma}}  \in \cA(O_{n, x, \gamma})\widehat\otimes \cK $$
is the restriction on a subset for all  $x, y\in \zdn, \nn, $ and $\gamma\in B_{X_n}(x, l_n)\cap \Gamma_n$.  On the other hand, it follows from the support condition (5) in Definition 5.5, there
exists $\bar r>0$ such that
$$\supp\big( T^{(n)}(x, y) \big) \subseteq B_{\mathbb{R}_+\times V_n} \big(  t_x(x)(s(x)) \,, \, \bar{r} \big) \,.$$
Hence, $T_\gamma^{(n)}(x, y) =0$ whenever
$$d\Big( t_x(\gamma)(s(\gamma)) \,, \, t_x(x)(s(x)) \Big) > \bar{r}+r$$
in $V_n$. It follows that there exists $S>0$ (since the section $s$ is locally a coarse embedding) such that $T_\gamma^{(n)}(x, y) =0$ whenever $d(x, \gamma)>S$ for all
$x, y\in \zdn, \gamma\in \Gamma_n,  \nn $. Now define
\begin{equation}\label{rrr}
S_\gamma^{(n)}(x, y)=\big( t_{\gamma  x} \big)_*\big( T_\gamma^{(n)}(x, y) \big) \,.
\end{equation}
Since $\big( T^{(0)}, \cdots, T^{(n)}, \cdots  \big)$ has finite propagation, there exists $N>0$ large enough such that, if $n\geq N$, then $S_\gamma^{(n)}(x, y)$ is well-defined, and
$$S_\gamma^{(n)}(x, y)\in \cA(O_{n, \gamma, \gamma})\widehat\otimes \cK$$
  for all $x, y\in \zdn\cap B_{\pdxn}(\gamma, S)$ with $\gamma\in \Gamma_n\cap B_{X_n}(x, l_n)$ and $n\geq N$.
Therefore, if we write $B(\gamma, S)=B_{\pdxn}(\gamma, S)$ and view a function as a matrix, we have
\[
\begin{array}{rcl}
S_\gamma^{(n)}  &  =  &  \Big[ S_\gamma^{(n)}(x, y)  \Big] _ {(x, y)\in (\zdn\times \zdn) \cap (B(\gamma, S)\times B(\gamma, S))}  \\
  & \in & \mathbb{C}\big[ B(\gamma, S)  \big]   \widehat\otimes \cA(O_{n, \gamma, \gamma})  \\
  & \subseteq   & C^*_{\max}\big( B(\gamma, S) \big)  \widehat\otimes \cA(O_{n, \gamma, \gamma}) \,.
\end{array}
\]
Define $S^{(n)}=0$ for all $n\leq N-1$ and define
$$S^{(n)}=\bigoplus_{\gamma\in \Gamma_n} S_\gamma^{(n)} \in  \bigoplus_{\gamma\in \Gamma_n} C^*_{\max}\big( B(\gamma, S) \big)  \widehat\otimes \cA(O_{n, \gamma, \gamma})$$
for $n\geq N$. Then the class $S=[( S^{(0)}, \cdots, S^{(n)}, \cdots )]$ is an element of
$$A_{\infty}\big[ \, \big( B(\gamma, S); \, \gamma\in \Gamma_n \big)_\nn \big] \, .$$
Now the correspondence $T\mapsto S$ extends to a $*$-isomorphism (see also the proof of Lemma 3.9 in \cite{SpWi} for essentially the same arguments which can be used to show that the norms in these two $C^*$-algebras agree), as desired in the first item. The remainder isomorphisms in this Proposition follow from the first one.
\hfill{$\Box$}
\par
Now we recall the notion of strong Lipschitz homotopy \cite{Yu97, Yu98, Yu00}.
\par
Let $(Y_\gamma)_{\ggnn}$ and $(\Delta_\gamma)_{\ggnn}$ be two families of uniformly bounded closed subsets of $\pdxn$, $\nn$,
for some $d\geq 0$, such that $\gamma\in Y_\gamma$, $\gamma\in \Delta_\gamma$ for all $\gamma\in \Gamma_n, \nn$.  A map
$$g: \; \bigsqcup_{{\gamma\in \Gamma_n, \nn}} Y_\gamma \longrightarrow \bigsqcup_{{\gamma\in \Gamma_n, \nn}} \Delta_\gamma $$
is said to be {\em Lipschitz} if
\\ \indent (1) $g(Y_\gamma)\subseteq \Delta_\gamma$ for each ${\gamma\in \Gamma_n, \nn}$;
\\ \indent (2) there exists a constant $c>0$, independent of $\ggnn$, such that
$$ d(g(x),g(y))\leq c \cdot d(x,y) $$
for all $x,y\in Y_\gamma$, $\ggnn$.
\par
Let $g_1,g_2$ be two Lipschitz maps from $\bigsqcup_{{\gamma\in \Gamma_n, \nn}} Y_\gamma$ to $\bigsqcup_{{\gamma\in \Gamma_n, \nn}} \Delta_\gamma$.
We say $g_1$ is {\it strongly Lipschitz homotopy} equivalent to $g_2$ if there exists a continuous map
$$ F: \; [0,1]\times \Big( \bigsqcup_{{\gamma\in \Gamma_n, \nn}} Y_\gamma \Big)\longrightarrow \bigsqcup_{{\gamma\in \Gamma_n, \nn}} \Delta_\gamma $$
such that
\\ \indent (1) $F(0,x)=g_1(x)$, $F(1,x)=g_2(x)$ for all $x\in\bigsqcup_{{\gamma\in \Gamma_n, \nn}} Y_\gamma$;
\\ \indent (2) there exists a constant $c>0$ for which $ d(F(t,x),F(t,y))\leq c \cdot d(x,y)$ for all $x,y\in Y_\gamma$, $\ggnn$, $t\in [0,1]$;
\\ \indent (3) $F$ is equicontinuous in $t$, i.e., for any $\varepsilon>0$ there exists $\delta>0$ such that $d(F(t_1,x),F(t_2,x))<\varepsilon $ for all
$x\in\bigsqcup_{{\gamma\in \Gamma_n, \nn}} Y_\gamma$ whenever $|t_1-t_2|<\delta$.
\par
We say $(Y_\gamma)_{\ggnn}$ is {\em strongly Lipschitz homotopy} equivalent to $(\Delta_\gamma)_{\ggnn}$ if there exist Lipschitz maps
$g_1:\bigsqcup Y_\gamma\to \bigsqcup \Delta_\gamma$ and $g_2:\bigsqcup \Delta_\gamma\to \bigsqcup Y_\gamma$ such that $g_1g_2$ and $g_2g_1$ are
respectively strongly Lipschitz homotopy equivalent to identity maps.
\par
Define $A^*_{L, 0, \infty}\big( (Y_\gamma; \, \gamma\in \Gamma_n)_\nn  \big)$ to be the $C^*$-subalgebra of
$$A^*_{L,  \infty}\big( (Y_\gamma; \, \gamma\in \Gamma_n)_\nn  \big)$$
consisting of those functions $f$ such that $f(0)=0$.
\par
{\bf Lemma 6.12.} (cf. \cite{Yu00})  {\it
If $(Y_\gamma)_{\ggnn}$ is strongly Lipschitz homotopy equivalent to $(\Delta_\gamma)_{\ggnn}$, then $K_*\big(A^*_{L, 0, \infty}\big( (Y_\gamma; \, \gamma\in \Gamma_n)_\nn  \big)\big)$ is isomorphic to $K_*\big(A^*_{L, 0, \infty}\big( (\Delta_\gamma; \, \gamma\in \Gamma_n)_\nn  \big)\big)$.
}
\hfill{$\Box$}
\par
Let $e$ be the evaluation homomorphism from $A^*_{L,  \infty}\big( (Y_\gamma; \, \gamma\in \Gamma_n)_\nn  \big)$   to $A^*_{\infty}\big( (Y_\gamma; \, \gamma\in \Gamma_n)_\nn  \big)$  given by
$e(f)=f(0)$.
\par
{\bf Proposition 6.13.} (cf. \cite{Yu00}) {\it For any $d\geq 0$, let $\Delta_\gamma$ be a simplex in $\pdxn$ for all $\ggnn$ with $\gamma\in \Delta_\gamma$. Then the evaluation map
$$ e_*:\, K_*\big(A^*_{L, \infty} \big( (\Delta_\gamma; \gamma\in \Gamma_n)_\nn \big)\big) \rightarrow  K_*\big(A^*_{\infty} \big( (\Delta_\gamma; \gamma\in \Gamma_n)_\nn \big)\big) $$
is an isomorphism.
}
\par
{\bf Proof.} (cf.  \cite{Yu00}) Note that $(\Delta_\gamma)_{\ggnn}$ is strongly Lipschitz homotopy equivalent to $(\{\gamma\})_{\ggnn}$. By an argument of Eilenberg swindle,  we have
$K_*\big(A^*_{L,0, \infty}\big( (\{\gamma\}; \gamma\in \Gamma_n )_\nn \big)\big)=0$. Consequently, Proposition 6.13 follows from Lemma 6.12 and the six-term exact sequence of $C^*$-algebra $K$-theory.
\hfill{$\Box$}
\par
We are now ready to give a proof to Theorem 6.7.
\par
{\bf Proof of Theorem 6.7. }  By Proposition 6.11, we have the following commutative diagram
\[
\xymatrix{
  \displaystyle\lim_{d\to\infty}    \culmpdxna_O         \ar[d]_{\cong}          \ar[r]^{\quad e}
            &    \displaystyle\lim_{d\to\infty}      \cumpdxna_O          \ar[d]^{\cong}                   \\
  \displaystyle\lim_{S\to\infty}     A^*_{L, \infty} \big( (\Delta_\gamma(S); \gamma\in \Gamma_n)_\nn \big)        \ar[r]^{\quad e}
            &    \displaystyle\lim_{s\to\infty}     A^*_{\infty} \big( (\Delta_\gamma(S); \gamma\in \Gamma_n)_\nn \big)
     }
\]
which induces the following commutative diagram on $K$-theory
\[
\xymatrix{
  \displaystyle\lim_{d\to\infty}    K_*\Big(  \culmpdxna_O  \Big)       \ar[d]_{\cong}          \ar[r]^{e_*}
            &    \displaystyle\lim_{d\to\infty}    K_*\Big(   \cumpdxna_O   \Big)       \ar[d]^{\cong}                   \\
  \displaystyle\lim_{S\to\infty}  K_*\Big(   A^*_{L, \infty} \big( (\Delta_\gamma(S); \gamma\in \Gamma_n)_\nn \big)   \Big)     \ar[r]^{e_*}
            &    \displaystyle\lim_{s\to\infty}    K_*\Big(    A^*_{\infty} \big( (\Delta_\gamma(S); \gamma\in \Gamma_n)_\nn \big)   \Big)
     }
\]
Now Theorem 6.7 follows from Proposition 6.13.
\hfill{$\Box$}
\subsection*{\it \S 6.3. Proof of Theorem 6.1}
We are now able to prove Theorem 6.1, the main result of this section, by gluing together the pieces we studied above.
\par
{\bf Proof of Theorem 6.1.}  For all $r>0$, define
$$O^{(r)}_{n, x} \;:=\bigcup_{\gamma\in B_{X_n}(x, l_n)} B_{\mathbb{R}_+\times V_n}\Big( t_x(\gamma)(s(\gamma)), \, r   \Big) $$
for all $x\in X_n, \nn$. Then
$$O^{(r)}:=\big( O^{(r)}_{n, x} \big)_{x\in X_n, \nn}$$
is a coherent system of open subsets. For any $d\geq 0$, if $r<r'$ then $O^{(r)}_{n, x}\subseteq O^{(r')}_{n, x}$ so that
$$\cumpdxna_{O^{(r)}} \subseteq \cumpdxna_{O^{(r')}}   \;,$$
$$\culmpdxna_{O^{(r)}} \subseteq \culmpdxna_{O^{(r')}} \;. $$
By the definition of these ideals, we have
$$\cumpdxna = \lim_{r\to \infty} \cumpdxna_{O^{(r)}}  \;, $$
$$\culmpdxna= \lim_{r\to \infty} \culmpdxna_{O^{(r)}} \;. $$
Consequently,  it suffices to show that,  for each $r_0> 0$, the evaluation map
\[
\begin{array}{ll}
e_*: &   \displaystyle \lim_{d\to \infty} K_*\Big( \displaystyle \lim_{r<r_0, r\to r_0} \culmpdxna_{O^{(r)}} \Big)             \\
&   \longrightarrow  \displaystyle \lim_{d\to \infty} K_*\Big(   \displaystyle \lim_{r<r_0, r\to r_0}  \cumpdxna_{O^{(r)}} \Big)
\end{array}
\]
is an isomorphism.
\par
Let $r_0>0$. Since $(X_n)_\nn$ has uniform bounded geometry, there exists $N>0$ such that $\# B(x, r_0)<N$ for all $x\in X_n, \nn$.
It follows that there exists an integer $J_{r_0}>0$ independent of $n$ such that we have the following decompositions for all $\nn$:
\begin{itemize}
\item $X_n=\bigcup_{j=1}^{J_{r_0}} \Gamma_n^{(j)}$ for $J_{r_0}$ subspaces $\Gamma_n^{(j)}$ of $X_n$;
\item $\Gamma_n^{(j)}\cap \Gamma_n^{(j')}=\emptyset$ whenever $j\not= j'$;
\item for any $x\in X_n$ and any distinct $\gamma, \gamma'\in \Gamma_n^{(j)}\cap B_{X_n}(x, l_n)$,
			$$d\Big( t_x(\gamma)(s(\gamma)), \, t_x(\gamma')(s(\gamma'))  \Big) >2r_0$$
		in $V_n$.
\end{itemize}
\par
For any $0<r<r_0$ and each $j\in \{1, 2, \cdots, J_{r_0}\}$, let
$$O^{(r, j)}_{n,x} = \bigcup_{\gamma\in \Gamma_n^{(j)}\cap B_{X_n}(x, l_n)}  B_{\mathbb{R}_+\times V_n}\big( t_x(\gamma)(s(\gamma)), \, r   \big) $$
for all $x\in X_n, \nn$. Then the system
$$O^{(r, j)}=\Big( O^{(r, j)}_{n,x} \Big)_{x\in X_n, \nn}$$
is coherent for each $j\in \{1, 2, \cdots, J_{r_0}\}$, and
$$O^{(r)}=\bigcup_{j=1}^{J_{r_0}} O^{(r, j)}.$$
For each $j\in \{1, 2, \cdots, J_{r_0}\}$,  let
$$\Gamma^{(j)}=\big( \Gamma_n^{(j)} \big)_\nn.$$
It is clear that the system $O^{(r, j)}$ or
$$O^{(r, j)} \cap \big( \cup_{i=1}^{j-1} O^{(r, i)}  \big)$$
is $(\Gamma^{(j)}, r)$-separate for each $j\in \{1, 2, \cdots, J_{r_0}\}$.
\par
By a construction of ``partition of unity'' as in the proof of Lemma 6.3 in \cite{Yu00}, for each $j\in \{1, 2, \cdots, J_{r_0}\}$, we have
$$\displaystyle \lim_{r<r_0, r\to r_0} \cumpdxna_{O^{(r, j)}} +\, \displaystyle \lim_{r<r_0, r\to r_0} \cumpdxna_{\cup_{i=1}^{j-1} O^{(r, i)}} $$
$$  =  \displaystyle \lim_{r<r_0, r\to r_0}  \cumpdxna_{\cup_{i=1}^{j} O^{(r, i)}} \;\; ,$$
$$ \displaystyle \lim_{r<r_0, r\to r_0} \cumpdxna_{O^{(r, j)}}\cap \, \displaystyle \lim_{r<r_0, r\to r_0} \cumpdxna_{\cup_{i=1}^{j-1} O^{(r, i)}} $$
$$ = \displaystyle \lim_{r<r_0, r\to r_0} \cumpdxna_{O^{(r, j)} \cap \big( \cup_{i=1}^{j-1} O^{(r, i)}  \big)} \;\; ,$$
$$\displaystyle \lim_{r<r_0, r\to r_0} \culmpdxna_{O^{(r, j)}}+\, \displaystyle \lim_{r<r_0, r\to r_0} \culmpdxna_{\cup_{i=1}^{j-1} O^{(r, i)}} $$
$$ = \displaystyle \lim_{r<r_0, r\to r_0}  \culmpdxna_{\cup_{i=1}^{j} O^{(r, i)}} \;\; ,$$
$$\displaystyle \lim_{r<r_0, r\to r_0} \culmpdxna_{O^{(r, j)}}\cap \,\displaystyle \lim_{r<r_0, r\to r_0} \culmpdxna_{\cup_{i=1}^{j-1} O^{(r, i)}} $$
$$  = \displaystyle \lim_{r<r_0, r\to r_0}  \culmpdxna_{O^{(r, j)} \cap \big( \cup_{i=1}^{j-1} O^{(r, i)}  \big)} \;\; .$$
\par
Now, Theorem 6.1 follows from Theorem 6.7, together with a Mayer-Vietoris sequence argument (cf. Sect. 3 of \cite{HRY}).
\hfill{$\Box$}

\section{Geometric analogue of Bott periodicity}
In this section, we shall define maps $\alpha$, $\beta$, $\alpha_L$, $\beta_L$ to build the following commutative diagram
\[
\xymatrix{
K_*\Big( \culmpdxn \Big) \ar[r]^{\quad e_*} \ar[d]^{(\beta_L)_*} & K_*\Big( \cumpdxn \Big) \ar[d]^{\beta_*} \\
K_*\Big( \culmpdxna \Big)  \ar@<1ex>[u]^{(\alpha_L)_*} \ar[r]^{\quad e_*}  & K_*\Big( \cumpdxna \Big )\;. \ar@<1ex>[u]^{\alpha_*}
}
\]
and show a geometric analogue of the Bott periodicity.
\par
The maps $\alpha$, $\beta$, $\alpha_L$, $\beta_L$ will be constructed as {\em asymptotic morphisms}. The reader is referred to
\cite{GHT, HK97, HK01, HKT, HG, Yu00} for backgrounds on asymptotic morphisms and the sources of most of the ideas behind the current section. In what follows, we shall use the graded formalism from \cite{HG}.

\subsection*{\it \S 7.1. The Bott maps $\beta$ and $\beta_L$}
In this subsection, we define the Bott maps  $\beta$ and $\beta_L$.
Recall that $(X_n)_\nn$ is a sequence of finite metric spaces with uniform bounded geometry which admit a fibred coarse embedding into Hilbert space as in Definition 5.4.
For each $d\geq 0$, $\zdn$ is a countable dense subset of $\pdxn$. And
$$V_n:= \mbox{ affine-span }\Big\{ t_x(z)(s(z)) \, \Big| \, z\in B(x, l_n), x\in X_n \Big\}$$
is a finite dimensional affine subspace of $H$ for each $\nn$.
\par
By Definition 5.1,  for each $x\in \zdn$, the inclusion of the $0$-dimensional affine subspace $\{t_x(x)(s(x))\}$ into $V_n$ induces a $*$-homomorphism
$$\beta(x): \, \cS\cong \cA\big( \{t_x(x)(s(x))\} \big)  \rightarrow \cA(V_n)$$
by the formula
$$\Big( \beta(x) \Big) (g)= g \big(   X\widehat\otimes 1+ 1\widehat\otimes C_{V_n,t_x(x)(s(x))}   \big)  \,,$$
where
$$C_{V_n,t_x(x)(s(x))}:\, V_n\rightarrow V_n^0 \subseteq \Cliff(V_n^0)$$
is the function $v\mapsto v- t_x(x)(s(x))\in V_n^0 \subseteq \Cliff(V_n^0)$ for all $v\in V_n$.
\par
{\bf Definition 7.1.} Let $d\geq 0$. For each $t\in [1, \infty)$, define a map
$$\beta_t:\, \cS\widehat\otimes \cupdxn \longrightarrow \cumpdxna$$
for $g\in \cS$, $\ttt\in \cupdxn$ by the formula
$$\beta_t\big( g\widehat\otimes T \big)=\Big[ \, \Big(      \big(  \beta_t\big( g\widehat\otimes T \big) \big)^{(0)}, \cdots, \big(  \beta_t\big( g\widehat\otimes T \big) \big)^{(n)}, \cdots \Big)\,\Big] \,,$$
where
$$\big(  \beta_t\big( g\widehat\otimes T \big) \big)^{(n)}(x, y)=\big( \beta(x) \big) (g_t) \widehat\otimes T^{(n)}(x, y)$$
for $x, y\in \zdn$, $\nn$, and $g_t(r)=g(\frac{r}{t})$ for all $r\in \mathbb{R}$.
\hfill{$\Box$}
\par
{\bf Definition 7.2.}  Let $d\geq 0$. For each $t\in [1, \infty)$, define a map
$$(\beta_L)_t:\, \cS\widehat\otimes \culpdxn \longrightarrow \culmpdxna$$
by the formula
$$\Big(  (\beta_L)_t (f) \Big) (\tau)= \beta_t\big(  f(\tau)  \big)$$
for $\tau\in \mathbb{R}_+=[0, \infty)$.
\hfill{$\Box$}
\par
{\bf Lemma 7.3.} {\it
For each $d\geq 0$, the maps $(\beta_t)_{t\geq 1}$ and $\big((\beta_L)_t\big)_{t\geq 1}$ extend to asymptotic morphisms
$$\beta:\; \cS\widehat\otimes \cumpdxn  \rightarrow \cumpdxna   \,,$$
$$\beta_L:\; \cS\widehat\otimes \culmpdxn  \rightarrow \culmpdxna  \,.$$
}
\par
{\bf Proof.} For any $x, y\in \zdn$ with $n\geq N$ for some $N>0$ large enough such that both  $x$ and $y$ are in the regions of trivializations of $t_x$ and $t_y$, define
$1$-dimensional affine subspaces of $V_n$:
$$W^x:=W_{\{x, y\}}(x):=\mbox{ affine-span } \Big\{ t_x(x)(s(x)),\;  t_x(y)(s(y)) \Big\} \,,$$
$$W^y:=W_{\{x, y\}}(y):=\mbox{ affine-span } \Big\{ t_y(x)(s(x)), \; t_y(y)(s(y)) \Big\} \,.$$
Then, by Definition 5.1 on $\beta_{ba}$, we have
$$\beta(x)=\beta_{V_n, W^x}\circ \beta_{W^x, \{ t_x(x)(s(x)) \}} \,,$$
$$\beta(y)=\beta_{V_n, W^y}\circ \beta_{W^y, \{ t_y(y)(s(y)) \}} \,.$$
Moreover, we have  $W^x=t_{xy}\big( W^y \big)$ and
$$\beta_{W^x, \{ t_x(y)(s(y)) \}}(g)=\big( t_{xy} \big)_* \Big(  \beta_{W^y, \{ t_y(y)(s(y)) \}}(g)  \Big)$$
for all $g\in \cS$, where recall that $t_{xy}=t_x\circ t_y^{-1}: H\to H$ and
$$\big( t_{xy} \big)_*: \;  \cA(W_C(y)) \rightarrow \cA(W_C(x))$$
mapping $g\widehat\otimes h$ to $g\widehat\otimes (t_{xy})_*(h)$, is defined in Definition 5.3 for all $C\subseteq B_{X_n}(x, l_n)\cap B_{X_n}(y, l_n)$.
\par
For the generators $g(x)=\frac{1}{x\pm i}$ of $\cS=C_0(\mathbb{R})$, it is standard argument to verify, for $t\in [1, \infty)$, that
\[
\begin{array}{ll}
  &  \Big\| \beta_{W^x, \{ t_x(x)(s(x)) \}}(g_t) - \big( t_{xy} \big)_* \Big(  \beta_{W^y, \{ t_y(y)(s(y)) \}}(g_t)    \Big)  \Big\|              \\
= &  \big\| \beta_{W^x, \{ t_x(x)(s(x)) \}}(g_t) -   \beta_{W^x, \{ t_x(y)(s(y)) \}}(g_t)      \big\|              \\
\leq & \frac{1}{t} \|  t_x(x)(s(x)) - t_x(y)(s(y)) \|     \\
\leq & \frac{1}{t} \rho_2(d(x, y)) \; \longrightarrow 0 \;\; (t\to \infty)  \,.
\end{array}
\]
It follows from an approximation argument, together with \cite{Yu00}(Lemma 7.6) and \cite{HG}(Lemma 3.2), that for all $d\geq 0$, $R>0, r>0, c>0, \varepsilon>0$, there exists
$t_0>1$ such that for all $x, y\in \zdn, \nn, $ with $d(x, y)\leq R$, for all $t\geq t_0$ and all $g\in \cS$ with
$$\supp(g)\subseteq [-r, r], \quad\quad\quad  \|g'\|_\infty\leq c  \,,$$
we have
$$\Big\| \big(  \beta(x) \big) (g_t) -\big( t_{xy} \big)_*\Big(   \big(  \beta(y) \big) (g_t)  \Big)   \Big\| < \varepsilon.$$
($(t_{xy})_*$ is as in Definition 5.5) As a result, the maps $\beta_t$ ($t\geq 1$) define a $*$-homomorphism $\beta$ from $\cS\widehat\otimes \cupdxn$ to the asymptotic $C^*$-algebra
$$\mathfrak{A}\Big(  \cumpdxna  \Big) := \frac
{C_b\Big( \big[ 1, \infty \big), \; \cumpdxna   \Big)}
{C_0\Big( \big[ 1, \infty \big), \; \cumpdxna   \Big)}
$$
satisfying $\|\beta(g\widehat\otimes T)\|\leq \|g\|\cdot \|T\|$ for all $g\in \cS$ and $T\in \cupdxn$. Hence, by the universality of the max norm and the max tensor product,
$\beta$ extends to a $C^*$-homomorphism from
$$\cS \widehat\otimes _{\max} \cumpdxn $$
to $\mathfrak{A}\big(  \cumpdxna  \big)$. Since $\cS$ is nuclear, we conclude that $\beta_t$ extends to an asymptotic morphism from $\cS \widehat\otimes \cumpdxn$ to
$\cumpdxna$. The case for $\beta_L$ is similar.
\hfill{$\Box$}

\subsection*{\it \S 7.2. The Dirac maps $\alpha$ and $\alpha_L$      }
In this subsection, we define the Dirac maps $\alpha$ and $\alpha_L$. To begin with, we recall the definition of  the Bott-Dirac operator on a finite dimensional
affine subspace of the Hilbert space $H$.
\par
Let $V$ be a finite dimensional affine subspace of $H$, and let $V^0$ be the corresponding finite dimensional linear space consisting of differences of elements of $V$. Let
$$L^2(V):=L^2(V, \Cliff(V^0))$$
be the graded infinite dimensional complex Hilbert space of square integrable $\Cliff(V^0)$-valued functions on $V$, where $V$ is endowed with the Lebesgue measure induced from the inner
product on $V^0\subseteq H$. The grading of $L^2(V)$ is inherited from the Clifford algebra $\Cliff(V^0)$, which is also considered as a graded Hilbert space in such a way that
the monomials $e_{i_1}\cdots e_{i_p}$ associated to an orthonormal basis of $V^0$ are orthonormal.
\par
The {\em  Dirac operator $D_V$ of $V$}  is the unbounded operator on $L^2(V)$ defined by the formula
$$\big( D_V\xi \big) (v) =\sum_{i=1}^n (-1)^{\mbox{degree} (\xi)} \frac{\partial \xi}{\partial x_i} (v) \cdot e_i$$
where $\{e_1, \cdots, e_n\}$ is an orthonormal basis for $V^0$, $\xi\in L^2(V)$, $v\in V$, and $\{x_1, \cdots, x_n\}$ is the coordinates dual to  $\{e_1, \cdots, e_n\}$.
The domain of $D_V$ is the Schwartz subspace of $L^2(V)$. Note that $D_V$ does not depend on the choice of an orthonormal basis for $V^0$ (and does not depend on a basis point in $V$).
\par
The {\em Clifford operator $C_{V, v_0}$ of $V$ at $v_0\in V$} is an unbounded operator on $L^2(V)$ defined by the formula
$$\big( C_{V, v_0} \xi \big) (v) = (v-v_0) \cdot \xi(v)$$
for all $v\in V$, $\xi\in L^2(V)$, where the multiplication $\cdot$ is the Clifford multiplication of $v-v_0\in V^0\subset \Cliff(V^0)$ and $\xi(v)\in \Cliff(V^0)$.
The domain of $C_{V,v_0}$ is also the Schwartz subspace of $L^2(V)$. Note that the Clifford operator  $C_{V, v_0}$  depends on  the choice of the base point $v_0$.
If $V$ is actually a linear subspace, and $v_0=0\in V$, then we also denote $C_V:=C_{V, 0}$.
\par
The {\em  Bott-Dirac operator of  $V$ at $v_0\in V$} is
$$B_{V, v_0} =D_V+C_{V, v_0}$$
with domain again the Schwartz subspace of $L^2(V)$. Denote by $\xi_{V, v_0}$ the unit vector of $L^2(V)=L^2(V, \Cliff(V^0))$:
$$\xi_{V, v_0}(v)= \pi^{-\frac{n}{4}}\cdot e^{-\frac{\|v-v_0\|^2}{2}} \; ,$$
 where $n=\dim(V^0)$. Then $\xi_{V, v_0}$ spans the kernel of $B_{V, v_0}$, which is a $1$-dimensional subspace of $L^2(V)$.
If $V$ is a linear subspace, and $v_0=0\in V$, then we also denote $B_V:=B_{V, 0}$.
\par
Let $V_a, V_b$ be finite dimensional affine subspaces of $H$ such that $V_a\subseteq V_b$ and $V_b=V_{ba}^0\oplus V_a$. Then
$$L^2(V_b)\cong L^2(V_{ba}^0)\widehat\otimes L^2(V_a)$$
and $L^2(V_a)$ is regarded as a subspace of $L^2(V_b)$ via the isometry from $L^2(V_a)$ to $L^2(V_b)$ given by
$$\xi \mapsto \xi_{V_{ba}^0, 0}\widehat\otimes \xi.$$
\par
Recall also that an affine isometry $t: V\to V$ canonically induces a unitary on $L^2(V)$ and moreover, a $*$-isomorphism
$$t_*: \cK(L^2(V))\rightarrow \cK(L^2(V))$$
by conjugation with the unitary.
\par
{\bf We now come back to the case of interest.} Let $(X_n)_\nn$ be a sequence of finite metric spaces with uniform bounded geometry which admit a fibred coarse embedding
into Hilbert space as described in Definition 5.4. Recall in Definition 5.5, for all $\nn$,
$$V_n=\mbox{ affine-span }\big\{ t_x(z)(s(z)) \,|\, z\in B_{X_n}(x, l_n), x\in X_n \big\}.$$
Define
$$E_n:= \mbox{ linear-span } \Big\{ V_n,\, 0 \Big\} \subseteq H  \,,$$
which is a  finite dimensional {\em linear} subspace of $H$ containing $V_n$.
Denote
$$L^2_n:=L^2(E_n, \Cliff(E_n))$$
and denote by $\cK(L^2_n)$ the graded $C^*$-algebra of all compact operators on the graded Hilbert space $L^2_n$ for all $\nn$.
\par
{\bf Definition 7.4.} For all $x, z\in X_n (\nn) $ with $B(x, l_n)\cap B(z, l_n)\not=\emptyset$, we define an isomorphism
$$(t_{xz})_*:\, \cK(L^2_n)\widehat\otimes \cK \rightarrow \cK(L^2_n)\widehat\otimes \cK$$
as follows. Denote
$$W_x:=\mbox{ affine-span }\big\{ t_x(w)(s(w)) \,\big| \, w\in  B(x, l_n)\cap B(z, l_n)  \,  \}, $$
$$W_z:=\mbox{ affine-span }\big\{ t_z(w)(s(w)) \,\big| \, w\in  B(x, l_n)\cap B(z, l_n)  \,  \}. $$
Then $W_x=t_{xz}(W_z)$. Denote $W_x^\perp=E_n\ominus W_x$ and $W_z^\perp=E_n\ominus W_z$ the linear orthogonal complements of $W_x, W_y$ in $E_n$. Choose a unitary operator
$U_{xz}: \, W_z^\perp \rightarrow W_x^\perp $. Then
$$U_{xz}\oplus t_{xz}: \, E_n=W_z^\perp\oplus W_z \longrightarrow E_n=W_x^\perp\oplus W_x$$
is an affine isometry from $E_n$ onto $E_n$. We define
$$(t_{xz})_*:=(U_{xz}\oplus t_{xz})_*\widehat\otimes 1: \; \cK(L^2_n)\widehat\otimes \cK \rightarrow \cK(L^2_n)\widehat\otimes \cK .$$
Moreover, in Rips complexes for each $d\geq 0$, if $x, z\in \zdn\subseteq \pdxn, \nn$, with $x\in \Star(\bar x)$ and $z\in \Star(\bar z)$, we define $(t_{xz})_*=(t_{\bar{x}\bar{z}})_*$.
\par
{\bf Remark 7.5.} Let  $\xi_{W_x^\perp}$ and $\xi_{W_z^\perp}$ be the unit vectors in the kernels of the Bott-Dirac operators at the origin of the linear spaces $L^2(W_x^\perp)$ and $L^2(W_z^\perp)$.
Then  $\xi_{W_x^\perp} = U_{xz}\big( \xi_{W_z^\perp} \big)$. This fact implies in the sequel actual applications,  $(t_{xz})_*$ does not depend on the choice of the unitary $U_{xz}$.
\par
{\bf Definition 7.6.} For each $d\geq 0$, define $\cupdxnk$ to be the set of all equivalence classes $\ttt$ of sequences
	$$(T^{(0)}, \cdots, T^{(n)}, \cdots)$$
described as follows:
\begin{itemize}
\item[(1)] $T^{(n)}$ is a bounded function from $\zdn \times \zdn$ to $\cK(L^2_n)\widehat\otimes\cK$ for all $\nn$;
\item[(2)] for any bounded subset $B\subset \pdxn$, the set
			$$\{ (x, y)\in B\times B \cap \zdn\times \zdn \; | \; T^{(n)} (x, y)\not= 0 \}$$
			is finite;
\item[(3)] there exists $L>0$ such that
			$$\#\{y\in \zdn \; | \; T^{(n)}(x, y)\not= 0 \}<L, \quad \#\{y\in \zdn \; | \; T^{(n)}(y, x)\not= 0 \}<L $$
			for all $x\in \zdn, \nn$;
\item[(4)] there exists $R>0$ such that $T^{(n)}(x, y)=0$ whenever $d(x, y)>R$ for $x, y\in \zdn, \; \nn$.
\end{itemize}
The equivalence relation $\sim$ on these sequences is defined by
$$(T^{(0)}, \cdots, T^{(n)}, \cdots) \sim (S^{(0)}, \cdots, S^{(n)}, \cdots)$$ if and only if
$$\lim_{n\to \infty} \sup_{x, y\in \zdn} \big\| T^{(n)}(x, y) - S^{(n)}(x, y) \big\|_\cK =0.$$
The product structure for $\cupdxnk$ is defined as follows. For any two elements $\ttt$ and $S=\big[ (S^{(0)}, \cdots, S^{(n)}, \cdots) \big]$ in $\cupdxnk$, their product is defined to be
$$\big[ \big( (TS)^{(0)}, \cdots, (TS)^{(n)}, \cdots \big) \big], $$
where there exists a sufficiently large $N\in \mathbb{N}$ depending on the propagation of the two representative sequences,
such that $ (TS)^{(n)}=0$ for all $n=0, 1, 2, \cdots, N-1$ and
$$(TS)^{(n)}(x, y)=\sum_{z\in \zdn} \Big( T^{(n)}(x, z) \Big) \cdot \Big( \big(t_{xz}\big)_* \big( S^{(n)}(z, y) \big) \Big)$$
for all $x, y\in \zdn, \; n\geq N$, where $\big(t_{xz}\big)_*$ is defined as in Definition 7.4.
\par
The $*$-structure for $\cupdxnk$ is defined by the formula
$$ \big[ \big( T^{(0)}, \cdots, T^{(n)}, \cdots \big) \big]^*=   \big[ \big( (T^*)^{(0)}, \cdots, (T^*)^{(n)}, \cdots \big) \big] \,,$$
where
$$(T^*)^{(n)}(x, y)=\big( t_{xy} \big)_* \Big( \big( T^{(n)}(y, x)\big)^*\Big)$$
for all but finitely many $n$, and $0$ otherwise.
\par
Now, $\cupdxnk$ is made into a $*$-algebra by using the additional usual matrix operations. Define $\cumpdxnk$ to be the completion of $\cupdxnk$ with respect to the norm
\begin{footnotesize}
$$\big\|T\big\|_{\max} :=\sup\Big\{ \big\|\phi(T)\big\|_{\cB(H_\phi)} \; \Big | \; \phi:\cupdxnk \to \cB(H_\phi), \mbox{ a $*$-representation} \Big\}. $$
\end{footnotesize}
\hfill{$\Box$}
\par
Similarly we define the localization algebra $\culmpdxnk$.
\par
{\bf Definition 7.7.} Let $d\geq 0$. For each $t\in [1, \infty)$ define a map
$$\alpha_t:\; \cupdxna\rightarrow \cumpdxnk$$
by the formula
$$\alpha_t(T)=\Big[ \Big(  \big(\alpha_t(T)\big)^{(0)}, \cdots,  \big(\alpha_t(T)\big)^{(n)}, \cdots     \Big)  \Big]$$
for $\ttt\in \cupdxna$, with
$$\big(\alpha_t(T)\big)^{(n)}(x, y)=\Big( \theta_t^k(x) \Big) \big( T_{1}^{(n)}(x, y)\big)$$
for all $x, y\in \zdn, \nn$, where
\begin{itemize}
\item[(1)] there exists $k\geq 0$ independent of $\nn$ as in condition (6) of Definition 5.5 such that
						$$ T^{(n)}(x, y) = \Big( \beta_{V_n, W_k(x)} \widehat\otimes 1 \Big) \big( T^{(n)}_1 (x, y) \big)$$
		for some $T_1^{(n)}(x, y)\in \cA(W_k(x))\widehat\otimes \cK$ of the form $\sum_{i=1}^K g_i\widehat\otimes h_i\widehat\otimes k_i$ where $g_i\in \cS$, $h_i\in \cC(W_k(x))$,
			$k_i\in \cK$ for $1 \leq i \leq K$.
\item[(2)] The map
		$$ \theta_t^k(x):\, \cA(W_k(x))\widehat\otimes \cK \rightarrow \cK(L^2_n)\widehat\otimes \cK$$
		is defined  by the formula
		$$\Big( \theta_t^k(x) \Big) \big( g\widehat\otimes h \widehat\otimes k \big)
			= g_t\big( B_{E_n\ominus W_k(x)} \widehat\otimes 1 + 1\widehat\otimes D_{W_k(x)}  \big)\big( 1\widehat\otimes M_{h_t} \big) \widehat\otimes k$$
		for all $g\in \cS$, $h\in \cC(W_k(x))$, $k\in \cK$, $t\geq 1$, $x\in \zdn, \nn$, where
	\begin{itemize}
	\item[$\bullet$]  $g_t(s)=g(s/t)$ for $s\in \mathbb{R}$.	
	\item[$\bullet$] $D_{W_k(x)}$ is the Dirac operator of the affine space
					$$W_k(x)=\mbox{ affine-span } \big\{ t_x(z)(s(z)) \, \big| \, z\in B_{\pdxn}(x, k) \big\} \subset E_n .$$
	\item[$\bullet$] $B_{E_n\ominus W_k(x)}$ is the Bott-Dirac operator at the origin $0$ of the linear orthogonal complement of $W_k(x)$ in $E_n$.
	\item[$\bullet$] for any $h\in \cC(W_k(x))$ and $t\geq 1$, the function $h_t\in \cC(W_k(x))$ is defined by
					$$ h_t(v)= h\Big( t_x(x)(s(x))+\frac{1}{t} \big( v-t_x(x)(s(x))   \big)  \Big)$$
				for all $v\in W_k(x)$.
	\item[$\bullet$] $M_{h_t}$ is the pointwise multiplication operator on $L^2(W_k(x), \Cliff(W_k(x)^0))$
					$$\Big( M_{h_t} \xi \Big) (v) = h_t(v)\cdot \xi(v)$$
				for all $\xi\in L^2(W_k(x))$ and $v\in W_k(x)$.
	\item[$\bullet$] $1\widehat\otimes M_{h_t}:\, L^2_n\to L^2_n$ is defined according to the tensor decomposition
					$$L^2_n\cong L^2(E_n\ominus W_k(x)) \widehat\otimes L^2(W_k(x)).$$
	\end{itemize}
\end{itemize}
\par
{\bf Remark 7.8.} (1) The fact that $\Big( \theta_t^k(x) \Big) \big( g\widehat\otimes h \widehat\otimes k \big) $ is in $\cK(L^2_n)\widehat\otimes \cK$ follows from the ellipticity of the Dirac operator
and the Rellich lemma, see (\cite{HKT} Definition 8). (2) By the proof of Proposition 4.2 in \cite{HKT}, we know that $\alpha_t(T)$ does not asymptotically depend on the choice of $k$.
\par
{\bf Definition 7.9.} For each $d\geq 0$ and $t\geq 1$, define
$$(\alpha_L)_t:\; \culpdxna\rightarrow \culmpdxnk$$
by the formula
$$\Big( (\alpha_L)_t(f) \Big)(\tau)=\alpha_t\big(f(\tau) \big)$$
for all $\tau\in [0, \infty)$.
\par
{\bf Lemma 7.10.} For each $d\geq 0$, the maps $(\alpha_t)_{t\geq 1}$ and $((\alpha_L)_t)_{t\geq 1}$ extend to asymptotic morphisms
$$\alpha:\; \cumpdxna\rightarrow \cumpdxnk \,,$$
$$\alpha_L:\; \culmpdxna\rightarrow \culmpdxnk  \,.$$
\par
{\bf Proof.} We first fix some notations. For a given large $\ell \geq 0$, and for any subset
$$C\subseteq B_{\pdxn}(x, \ell)\cap B_{\pdxn}(z, \ell)$$
where $x, z\in \zdn$  (for all but finitely many $n$)  with $d(x, z)<\ell$ , denote
$$W_C(x)=\mbox{ affine-span }\big\{ t_x(w)(s(w)) \, \big| \, w\in C \big\},$$
$$W_C(z)=\mbox{ affine-span }\big\{ t_z(w)(s(w)) \, \big| \, w\in C \big\}.$$
Then $W_C(x)=t_{xz}(W_C(z))$. For any affine subspace $W$ with $W_C(x)\subseteq W\subseteq E_n$, we identify $\cA(W_C(x))$ with a subalgebra of $\cA(W)$ via the map
$\beta_{W, W_C(x)}$ defined in Definition 5.1.
\par
{\it Step 1.}  For any $K>0$, $r>0$, $c>0$ and $x\in \zdn$, denote by
		$$\Big[  \cA(W_C(x))\widehat\otimes \cK \Big]_{K, r, c}$$
the subset of $\cA(W_C(x))\widehat\otimes \cK$ consisting of those elements of the form $\sum_{i=1}^K g_i\widehat\otimes h_i\widehat\otimes k_i$ where $g_i\in \cS$, $h_i\in \cC(W_C(x))$,
$k_i\in \cK$ such that (1) each $g_i$ is supported in $[-r, r]$; (2) each $g_i$ and $h_i$ is continuously differentiable and their derivatives satisfy $\|g'_i\|_\infty\leq c$ and
 $\|\nabla_v h_i\|\leq c$  for all $v\in W_C(x)$ such that $\|v-t_x(x)(s(x))\|\leq 1$.
\par
It follows from (\cite{Yu00} Lemma 7.5) and (\cite{HKT} Lemma 2.9) that for any $\varepsilon>0$ there exists $t_0>1$ such that for all $t\geq t_0$ and all
$a, b\in \Big[  \cA(W_C(x))\widehat\otimes \cK \Big]_{K, r, c} $, we have
$$\big\| \theta_t^\ell(x)(ab)-\theta_t^\ell(x)(a)\theta_t^\ell(x)(b) \big\|<\varepsilon, \quad\quad   \big\| \theta_t^\ell(x)(a^*)-\theta_t^\ell(x)(a)^* \big\|<\varepsilon. $$
\par
{\it Step 2.} The affine isometry $t_{xz}: W_C(z)\to W_C(x)$ induces a diagram
\[
\xymatrix{
\cA\big(W_C(z)\big)\widehat\otimes \cK  \ar[r]^{\quad \theta_t^\ell(z)} \ar[d]_{(t_{xz})_*}        &    \cK(L^2_n)\widehat\otimes \cK  \ar[d]^{(t_{xz})_*}             \\
\cA\big(W_C(x)\big)\widehat\otimes \cK  \ar[r]^{\quad \theta_t^\ell(x)}   &  \cK(L^2_n)\widehat\otimes \cK
}
\]
\par
It follows from the proof of (\cite{HKT} Proposition 4.2) and (\cite{Yu00} Lemma 7.2) that the diagram is asymptotically commutative in the sense that, for any
$R>0$ (relating propagations in $\cupdxna$) and any large $\ell$ (greater than $R+k$ where $k$ is as  Definition 5.5 (6)), for any  $K>0$, $r>0$, $c>0$ as above, and for any
$\varepsilon>0$, there exists $t_0>1$ such that for all $t\geq t_0$ and all $b\in \Big[  \cA(W_C(z))\widehat\otimes \cK \Big]_{K, r, c} $ , we have
$$\Big\|   \Big(\theta_t^\ell(x)\Big) \big( (t_{xz})_*(b) \big) -  (t_{xz})_* \Big(  \big(\theta_t^\ell(z)\big)(b)      \Big) \Big\|<\varepsilon,$$
whenever  $x, z\in \zdn$ satisfying  $d(x, z)\leq R$.
\par
{\it Step 3.} From the above facts it follows that the maps $(\alpha_t)_{t\geq 1}$ define a $*$-homomorphism  $\alpha$ from $\cupdxna$ to the asymptotic $C^*$-algebra
$$\mathfrak{A}\Big(  \cumpdxnk  \Big) := \frac
{C_b\Big( \big[ 1, \infty \big), \; \cumpdxnk   \Big)}
{C_0\Big( \big[ 1, \infty \big), \; \cumpdxnk   \Big)}
$$
By the universality of the max norm, $\alpha$ extends to a $*$-homomorphism from the $C^*$-algebra $\cumpdxna$ to $\mathfrak{A}$. The case for $\alpha_L$ is similar.
\hfill{$\Box$}

\subsection*{\it \S 7.3. A geometric analogue of the Bott periodicity in finite dimension}
Note that the asymptotic morphisms $\beta$, $\beta_L$, $\alpha$, $\alpha_L$  induce the following commutative diagram on $K$-theory:
\[
\xymatrix{
K_*\Big( \culmpdxn \Big) \ar[r]^{\quad e_*} \ar[d]^{(\beta_L)_*} & K_*\Big( \cumpdxn \Big) \ar[d]^{\beta_*} \\
K_*\Big( \culmpdxna \Big)  \ar[d]^{(\alpha_L)_*} \ar[r]^{\quad e_*}  & K_*\Big( \cumpdxna \Big )\;. \ar[d]^{\alpha_*}  \\
K_*\Big( \culmpdxnk \Big)  \ar[r]^{\quad e_*}  & K_*\Big( \cumpdxnk \Big )
}
\]
\par
In this subsection, we shall prove the following result, which is a geometric analogue of the Bott periodicity in finite dimension. The proof is adapted from (\cite{Yu00} Proposition 7.7).
\par
{\bf Theorem 7.11.} {\it
For each $d\geq 0$, the compositions $\alpha_*\circ \beta_*$ and $(\alpha_L)_*\circ (\beta_L)_*$ are the identity.
}
\par
{\bf Proof.} {\it Step 1.} Let $\gamma$ be the asymptotic morphism
$$\gamma:\, \cS\widehat\otimes \cumpdxn \rightarrow \cumpdxnk$$
defined by the formula
$$\gamma_t(g\widehat\otimes T)= \Big[ \, \Big(    \big(  \gamma_t(g\widehat\otimes T)   \big)^{(0)}, \cdots, \big(  \gamma_t(g\widehat\otimes T)   \big)^{(n)}, \cdots     \Big) \Big]$$
for all $g\in \cS$ and $\ttt\in \cupdxn$, where
$$\Big(  \gamma_t(g\widehat\otimes T)   \Big)^{(n)}(x, y) = g_{t^2}\big(  B_{E_n, t_x(x)(s(x))}  \big)\widehat\otimes T^{(n)}(x, y)$$
for all $t\geq 1$, $x, y\in \zdn, \nn$, and $B_{E_n, t_x(x)(s(x))}$ is the Bott-Dirac operator of $E_n$ at $t_x(x)(s(x))$.
\par
It follows from (\cite{HKT} Proposition 2.13 and Appendix) and (\cite{Yu00} Proposition 7.7) that the composition $\alpha\circ\beta$ is asymptotically equivalent to
$\gamma$ (see also \cite{SpWi} for a remark on compositions of asymptotic morphisms). Hence, $\alpha_*\circ\beta_*=\gamma_*$.
\par
{\it Step 2.} Let $\delta$ be the asymptotic morphism
$$\delta:\, \cS\widehat\otimes \cumpdxn \rightarrow \cumpdxnk$$
defined by the formula
$$\delta_t(g\widehat\otimes T)= \Big[ \, \Big(    \big(  \delta_t(g\widehat\otimes T)   \big)^{(0)}, \cdots, \big(  \delta_t(g\widehat\otimes T)   \big)^{(n)}, \cdots     \Big) \Big]$$
for all $g\in \cS$ and $\ttt\in \cupdxn$, where
$$\Big(  \delta_t(g\widehat\otimes T)   \Big)^{(n)}(x, y) = g_{t^2}\big(  B_{E_n, 0}  \big)\widehat\otimes T^{(n)}(x, y)$$
for all $t\geq 1$, $x, y\in \zdn, \nn$, and $B_{E_n, 0}$ is the Bott-Dirac operator of $E_n$ at the origin $0\in E_n$.
\par
For each $x\in \zdn, \nn$, let $U_x:\, L^2_n\to L^2_n$ be the unitary operator defined by
$$\Big( U_x\xi \Big)(v)=\xi\Big(  v-t_x(x)(s(x)) \Big)$$
for all $\xi\in L^2_n:=L^2(E_n, \Cliff(E_n))$ and $v\in E_n$. Then
$$U_x^{-1} \, B_{E_n, t_x(x)(s(x))} \, U_x \, = \,  B_{E_n, 0} \,\,.$$
\par
For each $s\in [0, 1]$, define an asymptotic morphism
$$\Phi^{(s)}:\, \cS\widehat\otimes \cumpdxn \rightarrow \cumpdxnk  \widehat\otimes \cM_2(\mathbb{C})          $$
by the formula
\[
\Big(  \Phi^{(s)}_t(g\widehat\otimes T)   \Big)^{(n)}(x, y)   =  \Big( U_x^{(s)}\Big)^{-1}
		\left[
		\begin{array}{cc}
			\Big(  \gamma_t(g\widehat\otimes T)   \Big)^{(n)}\big(x, y\big)     &      0                 \\
				0       &       0
		\end{array}
		\right]
U^{(s)}_x
\]
for all $t\geq 1$, $x, y\in \zdn, \nn$, where
\[
U^{(s)}_x =R(s) \left[
					\begin{array}{cc}
						U_x\widehat\otimes 1     &       0          \\
									0        &       1
					\end{array}
			\right]
		 R(s)^{-1},  \quad
R(s)= \left[
					\begin{array}{cc}
						\cos(\frac{\pi}{2}s)        &       \sin(\frac{\pi}{2}s)          \\
					     -\sin(\frac{\pi}{2}s)        &       \cos(\frac{\pi}{2}s)
					\end{array}
			\right].
\]
Then $\Phi^{(s)}$ is a homotopy between the asymptotic morphisms $\gamma=\Phi^{(1)}$ and $\delta=\Phi^{(0)}$. Hence, $\gamma_*=\delta_*$.
\par
{\it Step 3.} For each $\nn$, let $p^{(n)}$ be the projection of $E_n$ onto the 1-dimensional kernel of the Bott-Dirac operator $B_{E_n, 0}$ at the origin $0\in E_n$.
For each  $s\in [0, 1]$, define an asymptotic morphism
$$\Psi^{(s)}:\, \cS\widehat\otimes \cumpdxn \rightarrow \cumpdxnk   $$
by the formula
\[
\Big(  \Psi^{(s)}_t(g\widehat\otimes T)   \Big)^{(n)}(x, y) =\left\{
\begin{array}{ll}
g_{t^2}\big( \frac{1}{s} B_{E_n, 0}  \big) \widehat\otimes T^{(n)}(x, y),    &     \mbox{ if } s\in (0, 1];       \\
g(0)\cdot p^{(n)} \widehat\otimes T^{(n)}(x, y), &   \mbox{ if } s=0,
\end{array}
\right.
\]
for all $t\geq 1$, $x, y\in \zdn, \nn$.
Then $\Psi^{(s)}$ is a homotopy between the asymptotic morphism $\delta$ and the $*$-homomorphism
$$\sigma:\,  \cS\widehat\otimes \cumpdxn \rightarrow \cumpdxnk $$
defined by the formula
$$\sigma \Big( g\widehat\otimes \big[\big(  T^{(0)}, \cdots, T^{(n)}, \cdots    \big)\big] \Big)
	=g(0)\cdot  \Big[\Big(  p^{(0)}\widehat\otimes T^{(0)}, \cdots,  p^{(n)}\widehat\otimes T^{(n)}, \cdots    \Big)\Big] .$$
It is clear that $\sigma$ induces identity on $K$-theory. Therefore, we conclude that
$$\alpha_*\circ\beta_*=\gamma_*=\delta_*=\sigma_*=identity$$
on $K_*\Big(  \cumpdxn    \Big)$.
The case for $(\alpha_L)_*\circ(\beta_L)_*$ is similar.
\hfill{$\Box$}
\par
\vskip 5mm
\par
{\bf Summary of Proof of Theorem 1.1.}  It follows from Theorem 7.11 and  Theorem 6.1, together with an argument of diagram chasing, that
the evaluation map
\begin{small}
$$e_*: \lim_{d\to \infty} K_*\Big( \culmpdxn \Big) \longrightarrow \lim_{d\to \infty} K_*\Big( \cumpdxn \Big) $$
\end{small}
is an isomorphism for a sequence of finite metric spaces $(X_n)_\nn$ with uniform bounded geometry such that the coarse disjoint union $X=\bigsqcup_\nn X_n$ admits a fibred coarse embedding
into Hilbert space. That is, Theorem 4.4 holds, which implies Proposition 4.6 and, subsequently, Theorem 1.1.
\hfill{$\Box$}

\vskip 2cm

\begin{itemize}
\item[] Xiaoman Chen, \\ School of Mathematical Sciences, Fudan University, 220 Handan Road, Shanghai 200433, P. R. China.  \quad  E-mail: xchen@fudan.edu.cn

\item[] Qin Wang, \\ 1. Department of Applied Mathematics, Donghua University, 2999 Renmin North Road at Songjiang District, Shanghai, 201620, P. R. China.
\quad 2. Research Center for Operator Algebras, East China Normal University, 3663 Zhongshan North Road, 200062, Shanghai, P. R. China. 
\\ \quad  E-mail: qwang@dhu.edu.cn

\item[] Guoliang Yu \\  Department of Mathematics, Texas A\&M University, College Station, TX 77843-3368, USA.  \quad  E-mail: guoliangyu1@gmail.com

\end{itemize}

\end{document}